\documentclass[a4paper,12pt]{article}
\catcode`\"=13
\def"#1{\v#1}
\usepackage[dvips]{epsfig,graphicx}
\usepackage{amsmath,amsfonts,verbatim}

\newcommand{\CC}{{\mathbb C}}

\newcommand{\ths}{\theta^*}
\newcommand{\ov}{\overline}
\newcommand{\qed}{\hfill\hbox{\rule{3pt}{6pt}}}
\newcommand{\proof}{{\sc Proof. }}
\newcommand{\Label}{\label}
\newcommand{\G}{\Gamma}
\newcommand{\mat}{{\rm Mat}}

\newcommand{\As}{A^*}
\newcommand{\Ae}{A^{\varepsilon}}
\newcommand{\Es}{E^*}
\newcommand{\Ee}{E^{\varepsilon}}
\newcommand{\e}{\varepsilon}
\newcommand{\la}{\langle}
\newcommand{\ra}{\rangle}
\newcommand{\im}{\mathbf{i}}

\newcommand{\M}{{\rm Mat}_{d+1}(\CC)}
\newcommand{\MX}{\mat_X(\CC)}

\newtheorem{theorem}{Theorem}[section]
\newtheorem{lemma}[theorem]{Lemma}
\newtheorem{corollary}[theorem]{Corollary}

\newtheorem{definition}[theorem]{Definition}

\newtheorem{remark}[theorem]{Remark}
\newtheorem{notation}[theorem]{Notation}

\title{Leonard triples and hypercubes}

\author{"Stefko Miklavi"c\footnote{The author gratefully acknowledges support by the 
                                   US Department of State and the Fulbright Scholar Program, 
                                   and thanks the University of Wisconsin-Madison for hospitality 
                                   during his visit in Spring 2007.
                                   Supported in part by 
                                   ``Javna agencija za raziskovalno dejavnost Republike Slovenije'', program no. Z1-9614.}
 \\ Department of Mathematics and Computer Science \\
        Faculty of Education, University of Primorska \\ 6000 Koper, Slovenia \\
        miklavic@pef.upr.si}

\begin{document}
\maketitle

\begin{abstract}
  Let $V$ denote a vector space over $\CC$ with finite positive dimension.
  By a {\em Leonard triple} on $V$ we mean an ordered triple of linear operators on $V$ such that 
  for each of these operators there exists a basis of $V$ with respect to which the matrix 
  representing that operator is diagonal and the matrices representing the other two operators are 
  irreducible tridiagonal.
  
  Let $D$ denote a positive integer and let $Q_D$ denote the graph of the $D$-dimensional
  hypercube. Let $X$ denote the vertex set of $Q_D$ and let $A \in \MX$ denote the adjacency
  matrix of $Q_D$. Fix $x \in X$ and let $\As \in \MX$ denote the corresponding dual adjacency
  matrix. Let $T$ denote the subalgebra of $\MX$ generated by $A, \As$. We refer to $T$ as the
  {\em Terwilliger algebra of} $Q_D$ {\em with respect to} $x$. The 
  matrices $A$ and $A^*$ are related by the fact that $2 \im A = \As \Ae - \Ae \As$ and 
  $2 \im \As = \Ae A - A \Ae$, where $2 \im \Ae = A \As - \As A$ and $\im^2=-1$. 
  
  We show that the triple $A$, $\As$, $\Ae$ acts on each irreducible $T$-module as a 
  Leonard triple. We give a detailed description of these Leonard triples.
\end{abstract}

%%%%%%%%%%%%%%%%%%%%%%%%%%%%%%%%%%%%%%%%%%%%%%%%%%%%%%%%%%%%%%%%%%%%
\section{Introduction}
\label{sec:intro}

We start by recalling the definition of a Leonard pair. 
To describe this object we use the following terms.
Let $\CC$ denote the field of complex numbers.
A square matrix with entries in $\CC$ is called {\em tridiagonal} 
whenever each nonzero entry lies on either the diagonal, the subdiagonal, or the superdiagonal. 
A tridiagonal matrix is called {\em irreducible} whenever each entry on the subdiagonal is nonzero 
and each entry on the superdiagonal is nonzero.
\begin{definition}
\label{def:le_pair}
{\rm \cite[Definition 1.1]{ter5} Let $V$ denote a vector space over $\CC$ with finite positive 
dimension. By a {\em Leonard pair} on $V$ we mean an ordered pair of linear operators
$A: V \to V$ and $\As: V \to V$ which satisfy the conditions {\rm (i), (ii)} below.
\begin{itemize}
\item[{\rm (i)}]   There exists a basis for $V$ with respect to which the matrix representing $A$ is 
diagonal and the matrix representing $\As$ is irreducible tridiagonal.
\item[{\rm (ii)}]  There exists a basis for $V$ with respect to which the matrix representing $\As$ is 
diagonal and the matrix representing $A$ is irreducible tridiagonal.
\end{itemize}}
\end{definition}
Leonard pairs have been explored as linear algebraic objects, in connection with orthogonal 
polynomials, and as representations of certain algebras \cite{ter3}--\cite{ter6}.
The notion of a Leonard triple was introduced by Curtin in \cite{Cu}. We recall the definition.
\begin{definition}
\label{def:le_triple}
{\rm \cite[Definition 1.2]{Cu}
Let $V$ denote a vector space over $\CC$ with finite positive dimension. By a {\em Leonard triple} on 
$V$ we mean an ordered triple of linear operators $A: V \to V$, $\As: V \to V$, $\Ae:V \to V$ which 
satisfy the conditions {\rm (i)--(iii)} below.
\begin{itemize}
\item[{\rm (i)}]   There exists a basis for $V$ with respect to which the matrix representing $A$ is 
diagonal and the matrices representing $\As$ and $\Ae$ are each irreducible tridiagonal.
\item[{\rm (ii)}]  There exists a basis for $V$ with respect to which the matrix representing $\As$ is 
diagonal and the matrices representing $\Ae$ and $A$ are each irreducible tridiagonal.
\item[{\rm (iii)}] There exists a basis for $V$ with respect to which the matrix representing $\Ae$ is 
diagonal and the matrices representing $A$ and $\As$ are each irreducible tridiagonal.
\end{itemize}}
\end{definition}
Leonard triples are closely related to Leonard pairs.
Indeed, any ordered pair of distinct elements of a Leonard triple form a Leonard pair.
This allows us to take advantage of the literature concerning Leonard pairs in our
study of Leonard triples. 

\noindent
The isomorphism classes of Leonard pairs are in bijective correspondence with the polynomials in
the terminating branch of Askey-Wilson scheme \cite{ter9,ter8}. In particular, results concerning 
Leonard pairs also have interpretations as results concerning such polynomials. Consequently, 
results concerning Leonard triples also have interpretations as results concerning such polynomials.

\noindent
Leonard pairs play a role in representation theory \cite{ITT,Koe1,ter4,ter5,zhe} and
combinatorics \cite{Ca1,Cu3,CN,Go,ITT,ter3,ter1}. Consequently, also Leonard triples
play a role in representation theory and combinatorics.

\noindent
Leonard triples are also related to spin models \cite{Cu2}, generalized Markov problem in
number theory and the Poncelet problem in projective geometry \cite{KZ}.

\smallskip \noindent
In this paper we consider a situation in graph theory where Leonard triples arise naturally. The 
situation is described as follows. Let $D$ denote a positive integer, let $Q_D$ denote the graph of 
the hypercube with dimension $D$ (see Section \ref{sec:hyper} for formal definitions), and let $X$ 
denote the vertex set of $Q_D$. Let $\MX$ denote the $\CC$-algebra of matrices with entries in $\CC$ 
and with rows and columns indexed by $X$. Let $A \in \MX$ denote the adjacency matrix of $Q_D$. For 
the rest of this introduction fix $x \in X$. Let $\As=\As(x)$ denote the diagonal matrix in $\MX$ with 
$(y,y)$-entry $D-2 \partial(x,y)$ for $y \in X$, where $\partial$ denotes path-length distance. The 
matrix $\As$ is called the {\em dual adjacency matrix of} $Q_D$ {\em with respect to} $x$ \cite{ter3}. 
Let $T=T(x)$ denote the subalgebra of $\MX$ generated by $A, \As$. The algebra $T$ is known as the 
{\em Terwilliger algebra of} $Q_D$ {\em with respect to} $x$ \cite{ter3}. As we shall see, $A$ and 
$\As$ are related by the fact that $2 \im A = \As \Ae - \Ae \As$ and $2 \im \As = \Ae A - A \Ae$, where
$2 \im \Ae = A \As - \As A$ and $\im^2=-1$. We call $\Ae$ the {\em imaginary adjacency matrix of} 
$Q_D$ {\it with respect to} $x$. The matrices $A$, $\As$, $\Ae$ are similar; indeed we display an 
invertible matrix $P \in T$ such that $\As = P A P^{-1}$, $\Ae = P \As P^{-1}$, $A = P \Ae P^{-1}$. 

\noindent
Let $W$ denote an irreducible $T$-module. We show that the triple $A$, $\As$, $\Ae$ acts 
on $W$ as a Leonard triple. We give this triple a detailed description which is summarized as 
follows. Consider the three bases for $W$ afforded by Definition \ref{def:le_triple}.
For each of these bases we display two normalizations that we find attractive, and this yields 
six bases for $W$. We compute the matrices which represent $A$, $\As$, $\Ae$ with respect to these 
six bases. We display the inner products between each pair of these bases. We then display the 
transition matrices between each pair of these bases.
We remark that our paper extends the work of Go \cite{Go}.

%%%%%%%%%%%%%%%%%%%%%%%%%%%%%%%%%%%%%%%%%%%%%%%%%%%%%%%%%%%%%%%%%%%%%%%%%%%%%%%%%%%%%%%%%%%%%%%%%

\section{Preliminaries}
\label{sec:prelim}

In this section we review some definitions and basic results concerning distance-regular graphs. 
See the book of Brouwer, Cohen and Neumaier \cite{BCN} for more background information. 

\smallskip \noindent

Let $X$ denote a nonempty finite set. Let $\MX$ denote the $\CC$-algebra of matrices with entries in 
$\CC$ and with rows and columns indexed by $X$. For $B \in \MX$ let $B^t$ and $\ov{B}$ denote the 
transpose and the complex conjugate of $B$, respectively. Let $V=\CC^X$ denote the vector space over 
$\CC$ consisting of column vectors with entries in $\CC$ and rows indexed by $X$. We observe $\MX$ 
acts on $V$ by left multiplication. We refer to $V$ as the {\em standard module} of $\MX$. For 
$v \in V$ let $v^t$ and $\ov{v}$ denote the transpose and the complex conjugate of $v$, respectively. 
We endow $V$ with the Hermitean inner product $\la u,v \ra = u^t \overline{v} \; (u,v \in V)$. For 
$y \in X$ let $\hat{y}$ denote the vector in $V$ with a $1$ in the $y$ coordinate and $0$ in all other 
coordinates. Observe that $\{\hat{y} | y \in X\}$ is an orthogonal basis for $V$. The following will
be useful: for each $B \in \MX$ we have 
\begin{equation}
\label{cez}
\la u, B v \ra = \la \ov{B}^t u, v \ra \qquad \qquad (u,v \in V).
\end{equation}

Let $\G=(X,R)$ denote a finite, undirected, connected graph, without loops or multiple edges, with 
vertex set $X$, edge set $R$, path-length distance function $\partial$, and diameter 
$D:=\max \{\partial(x,y) |$ $\: x,y \in X\}$. For a vertex $x \in X$ and an integer $i \ge 0$ let
$\G_i(x)$ denote the set of vertices at distance $i$ from $x$. For an integer $k \ge 0$ we say 
$\Gamma$ is {\em regular with valency} $k$ whenever $|\G_1(x)|=k$ for all $x \in X$. We say $\G$ is 
{\em distance-regular} whenever for all integers $0 \le h,i,j \le D$ and all $x,y \in X$ with 
$\partial(x,y)=h$ the number
$$
  p_{ij}^h := | \G_i(x) \cap \G_j(y) |
$$
is independent of $x,y$. The constants $p_{ij}^h$ are known as the {\em intersection numbers} of $\G$. 
From now on we assume $\G$ is distance-regular with $D \ge 1$. For convenience set 
$c_i:=p_{1, i-1}^i \, (1 \le i \le D)$, $a_i:=p_{1i}^i \, (0 \le i \le D)$, 
$b_i:=p_{1, i+1}^i \, (0 \le i \le D-1)$,  $k_i:=p_{ii}^0 \, (0 \le i \le D)$, and $c_0=0$, $b_D=0$. 
We observe that $\G$ is regular with valency $k=k_1=b_0$ and that $c_i+a_i+b_i=k$ for $0 \le i \le D$.
By \cite[p. 127]{BCN} the following hold for $0 \le h,i,j \le D$: (i) $p_{ij}^h=0$ if one of $h,i,j$ is
greater than the sum of the other two; and (ii) $p_{ij}^h \ne 0$ if one of $h,i,j$ equals the sum of
the other two.

We now recall the Bose-Mesner algebra of $\G$. For $0 \le i \le D$ let $A_i$ denote the 
matrix in $\MX$ with entries
$$
  (A_i)_{x y} = \left\{ \begin{array}{ll}
                 1 & \hbox{if } \; \partial(x,y)=i, \\
                 0 & \hbox{if } \; \partial(x,y) \ne i \end{array} \right. \qquad (x,y \in X).
$$
We abbreviate $A=A_1$ and call this the {\em adjacency matrix of} $\G$. Let $M$ 
denote the subalgebra of $\MX$ generated by $A$. By \cite[p. 44]{BCN} the matrices
$A_0, A_1, \ldots, A_D$ form a basis for $M$. We call $M$ the 
{\em Bose-Mesner algebra of} $\G$. We observe that $M$ is commutative and semi-simple. 
By \cite[p. 45]{BCN} there exists a basis $E_0,E_1, \ldots, E_D$ for $M$ such that
\begin{eqnarray}
 E_0 = |X|^{-1} J, \label{pi1} \\
 E_0 + E_1 + \cdots + E_D = I, \label{pi2} \\
 E_i^t = E_i \;\; (0 \le i \le D), \label{pi3} \\
 \ov{E_i} = E_i \;\; (0 \le i \le D), \label{pi4} \\
 E_i E_j = \delta_{ij} E_i \;\; (0 \le i,j \le D), \label{pi5}
\end{eqnarray}
where $I$ and $J$ denote the identity and the all-ones matrix of $\MX$, respectively. For convenience 
we define $E_i=0$ if $i<0$ or $i>D$. The matrices $E_0,E_1, \ldots, E_D$ are known as the 
{\em primitive idempotents of} $\G$, and $E_0$ is called the {\em trivial} idempotent.
We recall the eigenvalues of $\G$. Since $E_0, E_1, \ldots, E_D$ is a basis for $M$, there exist 
scalars $\theta_0, \theta_1, \ldots, \theta_D \in \CC$ such that
\begin{equation}
\label{theta}
  A = \sum_{i=0}^D \theta_i E_i.
\end{equation}
Combining this with \eqref{pi1} and \eqref{pi5} we find $AE_i = E_iA =\theta_i E_i$ for 
$0 \le i \le D$ and $\theta_0=k$. The scalars $\theta_0, \theta_1, \ldots, \theta_D$ are real 
\cite[p. 197]{BI}. Observe that $\theta_0, \theta_1, \ldots, \theta_D$ are mutually distinct since $A$ 
generates $M$. We refer to $\theta_i$ as the {\em eigenvalue} of $\G$ {\em associated with} $E_i$. For 
$0 \le i \le D$ let $m_i$ denote the rank of $E_i$. We call $m_i$ the {\em multiplicity} of $\theta_i$. 

\smallskip \noindent
By \eqref{pi2}--\eqref{pi5}, 
\begin{equation}
\label{de}
  V = E_0V + E_1V + \cdots + E_DV \qquad \hbox{(orthogonal direct sum)}.
\end{equation}
By linear interpolation,
\begin{equation}
\Label{priden2}
E_i = \prod_{{0 \le j \le D \atop j \ne i}} {A - \theta_j I \over \theta_i-\theta_j} \qquad
(0 \le i \le D).
\end{equation}

\smallskip \noindent
We now recall the $Q$-polynomial property. Note that $A_i \circ A_j = \delta_{ij} A_i$ for
$0 \le i,j \le D$, where $\circ$ is the entry-wise multiplication. Therefore $M$ is closed under 
$\circ$. Thus there exist $q^h_{ij} \in \CC \; (0 \le h, i, j \le D)$ such that 
$$
  E_i \circ E_j = |X|^{-1} \sum_{h=0}^D q^h_{ij} E_h \qquad (0 \le i, j \le D).
$$
By \cite[Proposition 4.1.5]{BCN} the scalars $q_{ij}^h$ are real and nonnegative for
$0 \le h,i,j \le D$. The $q^h_{ij}$ are called the {\em Krein parameters} of $\G$. The graph $\G$ is 
said to be $Q$-{\em polynomial} (with respect to the given ordering $E_0, E_1, \ldots, E_D$ of the 
primitive idempotents) whenever the following hold for $0 \le h,i,j \le D$: (i) $q_{ij}^h=0$ if one of
$h,i,j$ is greater than the sum of the other two; and (ii) $q_{ij}^h \ne 0$ if one of $h,i,j$ equals
the sum of the other two. 

%%%%%%%%%%%%%%%%%%%%%%%%%%%%%%%%%%%%%%%%%%%%%%%%%%%%%%%%%%%%%%%%%%%%%%%%%%%%%%

\section{The Terwilliger algebra}
\label{sec:ter}

In this section we recall the dual Bose-Mesner algebra and the Terwilliger algebra of $\G$. 
For the rest of this section fix $x \in X$. 
For $0 \le i \le D$ let $\Es_i=\Es_i(x)$ denote the diagonal matrix in $\MX$ with entries
$$
  (\Es_i)_{y y} = \left\{ \begin{array}{lll}
                 1 & \hbox{if } \; \partial(x,y)=i, \\
                 0 & \hbox{if } \; \partial(x,y) \ne i \end{array} \right. \qquad (y \in X).
$$
We call $E_i^*$ the $i$th {\em dual idempotent of} $\G$ {\em with respect to} $x$. We observe
\begin{eqnarray}
  \Es_0+\Es_1 + \cdots + \Es_D = I, \Label{den1} \\
  E_i^{*t} = E_i^* \;\; (0 \le i \le D), \Label{den3} \\
  \ov{E_i^*}=E_i^* \;\; (0 \le i \le D), \Label{den2} \\
  E_i^* E_j^* =  \delta_{ij} E_i^*\;\; (0 \le i,j \le D). \Label{den4} 
\end{eqnarray}
By construction $E_0^*, E_1^*, \ldots, E_D^*$ are linearly independent. Let $M^*=M^*(x)$ denote the 
subalgebra of $\MX$ spanned by $E_0^*, E_1^*, \ldots, E_D^*$. We call $M^*$ the {\em dual Bose-Mesner 
algebra of} $\G$ {\em with respect to} $x$. We observe $M^*$ is commutative and semi-simple.

\smallskip \noindent
Assume $\G$ is $Q$-polynomial with respect to the ordering $E_0, E_1, \ldots, E_D$ of the 
primitive idempotents. Let $\As=\As(x)$ denote the diagonal matrix in $\MX$ with $(y,y)$-entry
\begin{equation}
\label{defas}
  \As_{yy}=|X| E_{xy} \quad (y \in X),
\end{equation}
where $E=E_1$. We call $\As$ the {\em dual adjacency matrix of} $\G$ {\em with respect to} $x$.
By \cite[Lemma 3.11]{ter3} $M^*$ is generated by $\As$. We recall the dual eigenvalues of $\G$.
Since $E_0^*, E_1^*, \ldots, E_D^*$ is a basis for $M^*$
there exist $\ths_0, \ths_1, \ldots, \ths_D \in \CC$ such that
\begin{equation}
\label{AsEE}
  \As = \sum_{i=0}^D \ths_i E_i^*.
\end{equation}
Combining this with \eqref{den4} we find $\As E_i^* = E_i^* \As =\ths_i E_i^* \; (0 \le i \le D)$. 
By \cite[Lemma 3.11]{ter3} $\ths_0, \ths_1, \ldots, \ths_D$ are real. The scalars 
$\theta^*_0,\theta^*_1,\ldots,\theta^*_D$ are mutually distinct since $A^*$ generates $M^*$.
Note that $\ths_i$ is an eigenvalue of $\As$ and $E_i^*V$ is the corresponding eigenspace 
$(0 \le i \le D)$. Using \eqref{den1}--\eqref{den4} we find 
\begin{equation}
\label{dual_direct_sum}
  V=E_0^*V + E_1^*V + \cdots + E_D^*V \qquad \hbox{(orthogonal direct sum)}.
\end{equation}
We call the sequence $\ths_0, \ths_1, \ldots, \ths_D$ the {\em dual eigenvalue sequence of } $\G$. 
Observe that for $0 \le i \le D$ the rank of $\Es_i$ is $k_i$. Therefore 
$k_i$ is the multiplicity with which $\ths_i$ appears as an eigenvalue of $\As$.  

\smallskip \noindent
By linear interpolation we obtain
\begin{equation}
\Label{dualpriden2}
E_i^* = \prod_{{0 \le j \le D \atop j \ne i}} {\As - \ths_j I \over \ths_i-\ths_j}.
\end{equation}
By \cite[Lemma 3.2]{ter3} the following hold for $0 \le h,j \le D$:
\begin{gather}
\label{ter:lem1i}
  E_j^* A E_h^* = 0 \; \hbox{ if and only if } \; p^h_{1j}=0; \\
\label{ter:lem1ii}
  E_j \As E_h = 0 \; \hbox{ if and only if } \; q^h_{1j}=0.
\end{gather}

\noindent
Let $T=T(x)$ denote the subalgebra of $\MX$ generated by $M$ and $M^*$. We call $T$ the 
{\em Terwilliger algebra of} $\G$ {\em with respect to} $x$ \cite[Definition 3.3]{ter3}.

\smallskip 
By a $T$-{\em module} we mean a subspace $W$ of $V$ such that $BW \subseteq W$ 
for all $B \in T$. Let $W$ denote a $T$-module. Then $W$ is said to be {\em irreducible} whenever
$W$ is nonzero and $W$ contains no $T$-modules other than $0$ and $W$. 

By construction $T$ is closed under the conjugate-transpose map so $T$ is semi-simple 
\cite[Lemma 3.4(i)]{ter3}. By \cite[Lemma 3.4(ii)]{ter3} $V$ decomposes into an orthogonal direct sum 
of irreducible $T$-modules. Let $W$ denote an irreducible $T$-module. By \cite[Lemma 3.4(iii)]{ter3} 
$W$ is the orthogonal direct sum of the nonvanishing $E_iW$ $(0 \le i \le D)$ and the orthogonal 
direct sum of the nonvanishing $E_i^*W \;(0 \le i \le D)$.
By the {\em endpoint} of $W$ we mean $\min \{i | 0 \le i\le D, \; \Es_i W \ne 0 \}$. By the 
{\em diameter} of $W$ we mean $|\{i | 0 \le i \le D, \; \Es_i W \ne 0\}| -1$. By the 
{\em dual endpoint} of $W$ we mean $\min \{i | 0 \le i\le D, \; E_i W \ne 0 \}$. By the 
{\em dual diameter} of $W$ we mean $|\{i | 0 \le i \le D, \; E_i W \ne 0\}| -1$. By 
\cite[Lemma 4.5]{ITT} the diameter and the dual diameter of $W$ coincide.
Let $r$ and $r^*$ denote the endpoint and the dual endpoint of $W$, respectively, and let $d$ denote 
the diameter of $W$. By \cite[Lemma 3.9(ii), Lemma 3.12(ii)]{ter3} the following hold for 
$0 \le i \le D$:
\begin{equation}
\label{ter:lem2i}
   E_i W \ne 0 \; \hbox{ if and only if } \; r^* \le i \le r^*+d,
\end{equation}
\begin{equation}
\label{ter:lem2ii}
  E_i^*W \ne 0 \; \hbox{ if and only if } \; r \le i \le r+d.
\end{equation}

Let $W$ denote an irreducible $T$-module. By \cite[Lemma 3.9, Lemma 3.12]{ter3} the following
are equivalent: (i) dim$(E_i W) \le 1$ for $0 \le i \le D$; (ii) dim$(\Es_i W) \le 1$ 
for $0 \le i \le D$. In this case $W$ is called {\em thin}.  

%%%%%%%%%%%%%%%%%%%%%%%%%%%%%%%%%%%%%%%%%%%%%%%%%%%%%%%%%%%%%%%%%%%%%%%%%%%%%%

\section{The hypercubes}
\label{sec:hyper}

In this section we recall the hypercube graph and some of its basic properties.
Let $D$ denote a positive integer, and let $\{0,1\}^D$ denote the set of sequences 
$(t_1, t_2, \ldots, t_D)$, where $t_i \in \{0,1\}$ for $1 \le i \le D$. Let $Q_D$ denote
the graph with vertex set $X=\{0,1\}^D$, and where two vertices are adjacent
if and only if they differ in exactly one coordinate. We call $Q_D$ the $D$-{\em cube} or a 
{\em hypercube}. The graph $Q_D$ is connected and for $y,z \in X$ the distance $\partial(y,z)$ 
is the number of coordinates at which $y$ and $z$ differ. In particular the diameter of $Q_D$ 
equals $D$. The graph $Q_D$ is bipartite with bipartition $X=X^+ \cup X^-$, where $X^+$ (resp.~$X^-$) 
is the set of vertices of $Q_D$ with an even (resp. odd) number of positive coordinates. By 
\cite[p.~261]{BCN} $Q_D$ is distance-regular with intersection numbers
\begin{equation}
\Label{int_num}
a_i=0, \qquad \qquad b_i=D-i, \qquad \qquad c_i=i, \qquad \qquad k_i={D \choose i}
\qquad \qquad (0 \le i \le D).
\end{equation}
Let $\theta_0 > \cdots > \theta_D$ denote the eigenvalues of $Q_D$. By \cite[p.~261]{BCN} these 
eigenvalues and their multiplicities are given by
\begin{equation}
\Label{eig_mult}
\theta_i = D-2i, \qquad \qquad \qquad \qquad m_i={D \choose i} \qquad \qquad \qquad \qquad
(0 \le i \le D).
\end{equation}
For $0 \le i \le D$ let $E_i$ denote the primitive idempotent of $Q_D$ associated with $\theta_i$. 
By \cite[Corollary 8.4.2]{BCN}, $Q_D$ is $Q$-polynomial 
with respect to $E_0, E_1, \ldots, E_D$. Moreover, it follows from 
\cite[Theorem 8.4.4]{BCN} that 
\begin{equation}
\label{enakost_pr_st}
p_{ij}^h = q_{ij}^h \qquad (0 \le h,i,j \le D).
\end{equation}
Since $Q_D$ is bipartite,
\begin{equation}
\label{odd}
p_{ij}^h=0 \; \hbox{ if } h+i+j \; \hbox{ is odd } \; (0 \le h,i,j \le D).
\end{equation}
Let $\ths_0, \ldots , \ths_D$ denote the dual eigenvalue sequence of $Q_D$ for the given 
$Q$-polynomial structure. Then $\ths_i=\theta_i \; (0 \le i \le D)$ \cite[Lemma 3.7]{Go}.
Fix $x \in X$. Let $\As=\As(x)$ denote the corresponding dual adjacency matrix,
and let $T=T(x)$ denote the corresponding Terwilliger algebra. By \eqref{AsEE} and since 
$\ths_i=D-2i$ we have
\begin{equation}
\label{AsEEvera}
\As_{yy} = D - 2 \partial(x,y) \qquad (y \in X).
\end{equation}
Since $Q_D$ is vertex-transitive the isomorphism class of $T(x)$ does not dependent on $x$. 
For notational convenience, for the rest of this paper we assume $x=(0,0, \ldots, 0)$.
By \cite[Theorem 4.2]{Go} we have
\begin{eqnarray}
A^{*2} A - 2 \As A \As + A A^{*2} = 4A, \Label{fundeq2} \\
A^2 \As - 2 A \As A + \As A^2 = 4\As. \Label{fundeq1} 
\end{eqnarray}

Let $W$ denote an irreducible $T$-module. By \cite[Theorem 6.3]{Go} $W$ is thin. 
By \cite[Theorem 6.3, Theorem 8.1]{Go} the endpoint and the dual endpoint of $W$ coincide.
Denoting this common value by $r$ we have $d=D-2r$ and $0 \le r \le D/2$,
where $d$ is the diameter of $W$ \cite[Theorem 6.3]{Go} .

%%%%%%%%%%%%%%%%%%%%%%%%%%%%%%%%%%%%%%%%%%%%%%%%%%%%%%%%%%%%%%%%%%%%%%%%%%%%%%%%%%%%%%%%%%%%
\section{The Cartesian product and the Kronecker product}
\label{kronecker}

In this section we recall the Cartesian product of graphs and the Kronecker product of matrices. 
For graphs $\G=(X,R)$ and $\G'=(X',R')$ let $\G \times \G'$ denote the graph with vertex set 
$X \times X'$, and with vertex $(u,u')$ being adjacent to vertex $(v,v')$ if and only if either $u=v$ 
and $u'$ is adjacent to $v'$ in $\G'$, or $u'=v'$ and $u$ is adjacent to $v$ in $\G$. We call 
$\G \times \G'$ the {\em Cartesian product} of $\G$ and $\G'$.

For $B \in \MX$ and $B' \in \mat_{X'}(\CC)$ let $B \otimes B'$ denote the matrix in 
$\mat_{X \times X'}(\CC)$, with a $((u,u'),(v,v'))$-entry equal to the $(u, v)$-entry of $B$ times the 
$(u', v')$-entry of $B'$. We call $B \otimes B'$ the {\em Kronecker product} of $B$ and $B'$.
By \cite[p.~107]{Eves}
\begin{equation}
\Label{kron1}
  (B_1 \otimes B'_1)(B_2 \otimes B'_2) = (B_1 B_2) \otimes (B'_1 B'_2).
\end{equation}
Also by \cite[p.~107]{Eves}
\begin{gather}
\Label{kron2}
  B \otimes (\gamma_1 B'_1 + \gamma_2 B'_2) = \gamma_1 B \otimes B'_1 + \gamma_2 B \otimes B'_2, \\
\label{kron2a}
(\gamma_1 B'_1 + \gamma_2 B'_2) \otimes B = \gamma_1 B'_1 \otimes B + \gamma_2 B'_2 \otimes B,
\end{gather}
where $\gamma_1, \gamma_2 \in \CC$. It is also known that
\begin{equation}
\Label{kron4}
 (B \otimes B')^t = B^t \otimes B'^t, \qquad \ov{B \otimes B'} = \ov{B} \otimes \ov{B'}.
\end{equation}
For a matrix $B \in \MX$ and an integer $r \ge 0$ let $B^{\otimes r}$ denote
$B \otimes B \otimes \cdots \otimes B$ ($r$ copies). We interpret $B^{\otimes 0}=1$.
Let $\G$ and $\G'$ be graphs with adjacency matrices $A$ and $A'$, respectively. By construction the 
adjacency matrix of $\G \times \G'$ is equal to $A \otimes I' + I \otimes A',$ where $I$ and $I'$ are 
the identity matrices of appropriate dimensions (see, for example, \cite[Section 12.4]{God}).
The hypercube $Q_D$ can be viewed as a Cartesian product $Q_1 \times Q_1 \times \cdots \times Q_1$ 
($D$ copies). By a simple induction argument we find that the adjacency matrix $A$ of $Q_D$ satisfies
\begin{equation}
\label{pr1}
A = \sum_{i=0}^{D-1} I_1^{\otimes i} \otimes A_1 \otimes I_1^{\otimes (D-1-i)},
\end{equation}
where $A_1$ and $I_1$ denote the adjacency matrix of $Q_1$ and the $2 \times 2$ identity matrix, 
respectively. The reader is invited to verify that a similar equation holds
for the dual adjacency matrix $\As$ of the hypercube $Q_D$:
\begin{equation}
\label{pr2}
\As = \sum_{i=0}^{D-1} I_1^{\otimes i} \otimes \As_1 \otimes I_1^{\otimes (D-1-i)}.
\end{equation}
%%%%%%%%%%%%%%%%%%%%%%%%%%%%%%%%%%%%%%%%%%%%%%%%%%%%%%%%%%%%%%%%%%%%%%%%%%%%%%

\section{The imaginary adjacency matrix of $\boldsymbol{Q_D}$}
\label{sec:fund}

In this section we define the imaginary adjacency matrix of the hypercube $Q_D$. 
We use the following notation.

\begin{notation}
\Label{epsilon}
{\rm Let $D$ denote a positive integer and let $Q_D=(X,R)$ denote the $D$-cube. Let $A$ denote the 
adjacency matrix of $Q_D$ and let $I$ denote the identity matrix in $\MX$. Fix 
$x = (0,0, \ldots, 0) \in X$, and let $\As=\As(x)$ and $T=T(x)$ denote the corresponding dual 
adjacency matrix and the Terwilliger algebra, respectively. We define the matrix $\Ae=\Ae(x)$ by
\begin{equation}
\label{Ae}
  \Ae = -\im (A \As - \As A) / 2,
\end{equation}
where $\im^2=-1$.}
\end{notation}
Note that $\Ae \in T$. We have the following observation.

\begin{lemma}
\label{fund:lem0}
With reference to Notation \ref{epsilon}, for $y,z \in X$ the (y,z)-entry of $\Ae$ is given by
$$
  \Ae_{yz} = \im (\partial(x,z)-\partial(x,y)) A_{yz}.
$$
\end{lemma}
\proof
Since $\As$ is diagonal we obtain $(A \As)_{yz} = A_{yz} \As_{zz}$ and 
$(\As A)_{yz} = \As_{yy} A_{yz}$. By \eqref{AsEEvera} $\As_{yy} = D-2 \partial(x,y)$ and 
$\As_{zz} = D-2 \partial(x,z)$. The result now follows using \eqref{Ae}. \qed

\medskip \noindent
Motivated by Lemma \ref{fund:lem0} we call $\Ae$ the {\em imaginary adjacency matrix of} $Q_D$
{\em with respect to} $x$. The following result will be useful.

\begin{lemma} 
\Label{fund:lem1}
With reference to Notation \ref{epsilon} the following {\rm (i)--(iii)} hold.
\begin{itemize}
\item[{\rm (i)}]   $A \As - \As A = 2\im \Ae$,
\item[{\rm (ii)}]  $\As \Ae - \Ae \As = 2\im A$,
\item[{\rm (iii)}] $\Ae A - A \Ae = 2\im \As$.
\end{itemize}
\end{lemma}
\proof (i) Immediate from \eqref{Ae}.

\smallskip 
\noindent
(ii), (iii) 
Eliminate $\Ae$ using \eqref{Ae} and simplify using \eqref{fundeq2} and 
\eqref{fundeq1}. \qed

\begin{lemma} 
\Label{fund:lem1a}
With reference to Notation \ref{epsilon} we have
\begin{equation}
\label{Y}
 \Ae = \sum_{i=0}^{D-1} I_1^{\otimes i} \otimes \Ae_1 \otimes I_1^{\otimes (D-1-i)},
\end{equation}
where $\Ae_1$ denotes the imaginary adjacency matrix of $Q_1$ and $I_1$ denotes the $2 \times 2$
identity matrix. 
\end{lemma}
\proof
Evaluate the left-hand side using \eqref{Ae} and then \eqref{pr1}, \eqref{pr2}. Simplify the
result using  \eqref{kron1}--\eqref{kron2a} and $\Ae_1 = -\im (A_1 \As_1 - \As_1 A_1)/2$. \qed
	        
%%%%%%%%%%%%%%%%%%%%%%%%%%%%%%%%%%%%%%%%%%%%%%%%%%%%%%%%%%%%%%%%%%%%%%%%%%%%%%
\section{The eigenvalues of the imaginary adjacency matrix}
\label{sec:eig}

In this section we describe the eigenvalues for the imaginary adjacency matrix of $Q_D$.
We begin with a definition.

\begin{definition}
\label{P}
{\rm With reference to Notation \ref{epsilon} we define
\begin{equation}
\label{PP}
  P=P_1^{\otimes D}
\end{equation}
where $P_1$ is the matrix with
rows and columns indexed by the set $\{0,1\}$, and with $(0,0)$-entry $1$, $(0,1)$-entry 1,
$(1,0)$-entry $-\im$, and $(1,1)$-entry $\im$.}
\end{definition}
\begin{lemma}
\label{PinT}
With reference to Notation \ref{epsilon} and Definition \ref{P} we have $P \in T$.
\end{lemma}
\proof 
Observe that the symmetric group $S_D$ acts on $X$ as a group of automorphisms of $Q_D$ by the rule
\begin{equation}
\label{action}
  (t_1, t_2, \ldots, t_D)^{\sigma} = (t_{\sigma(1)}, t_{\sigma(2)}, \ldots, t_{\sigma(D)})
  \qquad (\sigma \in S_D),
\end{equation}
where $(t_1, t_2, \ldots, t_D)$ is a vertex of $Q_D$. In fact $S_D$ is isomorphic to the stabilizer of the vertex $x$ in the full automorphism group of $Q_D$ \cite[Theorem 9.2.1]{BCN}. 

Observe that the above $S_D$-action on $X$ induces an $S_D$-action on $V$. For $\sigma \in S_D$ let 
$M_{\sigma}$ denote the matrix in $\MX$ that represents $\sigma$ with respect to the basis 
$\{ \hat{y} | y \in X \}$. By \cite[Subsection I.C]{Sch} $T$ is the centralizer algebra for $S_D$ on $V$, i.e.  
$$
  T=\{ M \in \MX | M M_{\sigma}=M_{\sigma} M \; \; \forall \sigma \in S_D\}.
$$
Pick distinct $1 \le i,j \le D$ and consider the involution $\sigma=(i,j) \in S_D$.
Since $S_D$ is generated by the involutions, to show $P \in T$ it suffices to show 
$P M_{\sigma}=M_{\sigma} P$. For $y,z \in X$ the $(y,z)$-entry of $M_\sigma$ is $1$ if 
$y^\sigma=z$ and $0$ if $y^\sigma \ne z$. By this and matrix multiplication the $(y,z)$-entry
of $P M_\sigma$ is $P_{y z^\sigma}$ and the $(y,z)$-entry of $M_\sigma P$ is $P_{y^\sigma z}$.
By \eqref{PP} and the definition of the Kronecker product $P_{y z^\sigma} = P_{y^\sigma z}$.
By these comments $(P M_\sigma)_{yz}=(M_\sigma P)_{yz}$. Therefore $P M_\sigma=M_\sigma P$ so 
$P \in T$. \qed

\medskip \noindent
We have an observation.

\begin{lemma}
\Label{fund:lem3a}
With reference to Notation \ref{epsilon} and Definition \ref{P} the following {\rm (i)--(iii)} hold.
\begin{itemize}
\item[{\rm (i)}]   $P \ov{P}^t = \ov{P}^t P = 2^D I$; 
\item[{\rm (ii)}]  $P^3 = 2^D(1-\im)^D I$;
\item[{\rm (iii)}] $P^{-1}$ exists.
\end{itemize}
\end{lemma}
\proof (i) We first observe that $P_1 \ov{P_1}^t = 2 I_1$, where $I_1$ denotes the $2 \times 2$
identity matrix. Using \eqref{kron1} and \eqref{kron4} we now obtain
$$
  P \, \ov{P}^t = P_1^{\otimes D} \big( \ov{P_1}^t \big)^{\otimes D} = 
  \big( P_1 \, \ov{P_1}^t \big)^{\otimes D} = (2 I_1)^{\otimes D} = 2^D I.
$$
(ii) We first observe that $P_1^3 = 2(1-\im)I_1$. Using \eqref{kron1} we now obtain
$$
  P^3 = P_1^{\otimes D} P_1^{\otimes D} P_1^{\otimes D} = \big( P_1^3 \big)^{\otimes D} =
  \big( 2 (1-\im) I_1 \big)^{\otimes D} = 2^D (1-\im)^D I.
$$
(iii) Clear from (ii) above. \qed

\begin{theorem}
\Label{fund:thm3a}
With reference to Notation \ref{epsilon} and Definition \ref{P},
$$
 \As = P A P^{-1}, \qquad \qquad \Ae = P \As P^{-1}, \qquad \qquad A = P \Ae P^{-1}.
$$
\end{theorem}
\proof 
To verify the equation on the left, evaluate $\As$ using \eqref{pr2} and $P A P^{-1}$ using
\eqref{kron1}, \eqref{pr1}, \eqref{PP}. Comparing the results using $\As_1 = P_1 A_1 P_1^{-1}$ we find 
$\As = P A P^{-1}$. The other two equations are similarly obtained. \qed

\begin{corollary}
\Label{fund:thm3}
With reference to Notation \ref{epsilon} the following {\rm (i)--(iii)} hold.
\begin{itemize}
\item[{\rm (i)}]   
The matrix $\Ae$ is diagonalizable.
\item[{\rm (ii)}]  
The eigenvalues of $\Ae$ are $D-2i \; (0 \le i \le D)$. 
\item[{\rm (iii)}] 
For $0 \le i \le D$ the eigenvalue $D-2i$ of $\Ae$ has multiplicity ${D \choose i}$.
\end{itemize}
\end{corollary}
\proof 
Use \eqref{eig_mult} and the equation on the right in Theorem \ref{fund:thm3a}. \qed

%%%%%%%%%%%%%%%%%%%%%%%%%%%%%%%%%%%%%%%%%%%%%%%%%%%%%%%%%%%%%%%%%%%
\section{The primitive idempotents of the imaginary adjacency matrix}
\label{sec:im_id}

In this section we consider the primitive idempotents for the imaginary adjacency matrix of $Q_D$.
We start with a definition.

\begin{definition}
\label{def:imid}
{\rm With reference to Notation \ref{epsilon}, for $0 \le i \le D$ we define $\Ee_i = P^{-1} E_i P$,
where $P$ is from Definition \ref{P} and $E_i$ is the $i$th primitive idempotent of $Q_D$.}
\end{definition}
Adopt Notation \ref{epsilon}. Using Definition \ref{def:imid}, \eqref{pi2}--\eqref{theta}, 
\eqref{priden2}, Lemma \ref{fund:lem3a}(i), and the equation on the right in Theorem \ref{fund:thm3a} 
we routinely find 
\begin{eqnarray}
\Ee_0+\Ee_1+ \cdots + \Ee_D = I, \Label{ien2} \\
\ov{\Ee_i}^t = \Ee_i \;\; (0 \le i \le D), \Label{ien3} \\
\Ee_i \Ee_j = \delta_{ij} \Ee_i \; (0 \le i,j \le D), \label{een4}\\
\Ae = \sum_{i=0}^D (D-2i) E_i^{\e}, \label{ien0} \\
\Ae \Ee_i = \Ee_i \Ae = (D-2i) \Ee_i \;\; (0 \le i \le D), \Label{ien1} \\
\Ee_i = \prod_{{0 \le j \le D \atop j \ne i}} {\Ae - \theta_j I \over \theta_i - \theta_j} \;
(0 \le i \le D). \label{imaginarypriden2} 
\end{eqnarray}
For $0 \le i \le D$ we note that $\Ee_i$ is the primitive idempotent of $\Ae$ associated with the 
eigenvalue $D-2i$. It follows from \eqref{imaginarypriden2} that $\Ee_i \in T$. 
We call $\Ee_i$ the $i$th {\em imaginary idempotent of} $Q_D$ {\em with respect to} $x$. For the rest 
of this paper we consider the following situation.
\begin{notation}
\label{blank1}
{\rm Let $D$ denote a positive integer and let $Q_D=(X,R)$ denote the $D$-cube. Let $A$ denote the
adjacency matrix of $Q_D$ and let $I$ denote the identity matrix in $\MX$. Let $V$ denote the 
standard module of $\MX$. Fix $x=(0,0, \ldots, 0) \in X$, and let $\As=\As(x)$ and $\Ae=\Ae(x)$ denote 
the corresponding dual adjacency and the imaginary adjacency matrix, respectively. Let $T=T(x)$ denote 
the corresponding Terwilliger algebra. Let $E_i$, $E^*_i$, $E^{\e}_i \: (0 \le i \le D)$ 
denote the primitive idempotents, the dual idempotents and the imaginary idempotents of $Q_D$, 
respectively. Let the matrix $P$ be as in Definition \ref{P}.}
\end{notation}
We record the following for later use.

\begin{lemma}
\label{ii:lem0}
With reference to Notation \ref{blank1} the following {\rm (i), (ii)} hold.
\begin{itemize}
\item[{\rm (i)}]  For $0 \le i \le D$, $\Ee_i V$ is the eigenspace of $\Ae$ for the
eigenvalue $D-2i$.
\item[{\rm (ii)}] $V = \Ee_0 V + \Ee_1 V + \cdots + \Ee_D V \qquad \hbox{(orthogonal direct sum)}$.
\end{itemize}
\end{lemma}
\proof (i) This follows from \eqref{ien1}. \\
(ii) Evaluate $V=P^{-1} I P V$ using \eqref{pi2} and Definition \ref{def:imid} to obtain
$V = \sum_{i=0}^D \Ee_i V$ (direct sum). This sum is orthogonal by \eqref{ien3} and \eqref{een4}.
\qed

\begin{lemma}
\Label{fund:lem4a}
With reference to Notation \ref{blank1} the following holds for $0 \le i \le D$:
$$
  \Es_i = P E_i P^{-1}, \qquad \qquad \Ee_i = P \Es_i P^{-1}, \qquad \qquad E_i = P \Ee_i P^{-1}.
$$
\end{lemma}
\proof 
To verify the left equation, simplify the right-hand side using \eqref{priden2} and the equation on 
the left in Theorem \ref{fund:thm3a}. Compare the result with \eqref{dualpriden2} and recall
$\theta_i = \ths_i$ for $0 \le i \le D$.
The other two equations are similarly obtained. \qed

\begin{corollary}
\Label{fund:cor4a}
With reference to Notation \ref{blank1}, let $W$ denote an irreducible $T$-module. 
Then the following holds for $0 \le i \le D$:
$$
  P E_i W = E_i^* W, \qquad P E_i^* W = E_i^{\e} W, \qquad P E_i^{\e} W = E_i W.
$$
\end{corollary}
\proof Immediate from Lemma \ref{fund:lem4a} and since $P^{-1} W = W$ by Lemma \ref{PinT}. \qed

\begin{corollary}
\Label{fund:cor4aa}
With reference to Notation \ref{blank1}, let $W$ denote an irreducible $T$-module. 
Then $\hbox{dim}(\Ee_i W) \le 1$ for $0 \le i \le D$.
\end{corollary}
\proof This follows from Corollary \ref{fund:cor4a} and since 
dim$(\Es_i W) \le 1$ for $0 \le i \le D$. \qed

\begin{lemma}
\label{im_end}
With reference to Notation \ref{blank1}, let $W$ denote an irreducible $T$-module
with endpoint $r$ and diameter $d=D-2r$. Then the following {\rm (i), (ii)} hold.
\begin{itemize}
\item[{\rm (i)}]
$\Ee_i W \ne 0$ if and only if $r \le i \le r+d$ $(0 \le i \le D)$;
\item[{\rm (ii)}]
$W = \Ee_r W + \Ee_{r+1} W + \cdots + \Ee_{r+d} W$ \hspace{1cm} (orthogonal direct sum).
\end{itemize}
\end{lemma}
\proof (i) By Corollary \ref{fund:cor4a} we have $\Ee_i W = P \Es_i W$. By \eqref{ter:lem2ii} we have 
$\Es_i W \ne 0$ if and only if $r \le i \le r+d$. The result follows.

\noindent
(ii)
Recall that $W=\sum_{i=0}^d \Es_{r+i} W$ (direct sum) by \eqref{ter:lem2ii}. Therefore 
$PW=\sum_{i=0}^d P \Es_{r+i} W$ (direct sum). Simplify this equation using $W=PW$ and 
Corollary \ref{fund:cor4a} to obtain $W = \sum_{i=0}^d \Ee_{r+i} W$ (direct sum). This sum is orthogonal by \eqref{ien3} and \eqref{een4}. \qed

\begin{lemma}
\Label{fund:lem4}
With reference to Notation \ref{blank1} the following {\rm (i)--(v)} are equivalent for
$0 \le h,j \le D$:
\begin{itemize}
\item[{\rm (i)}]   $p^h_{1j}=0$;
\item[{\rm (ii)}]  $E_h \Ae E_j = 0$; 
\item[{\rm (iii)}] $E_h^* \Ae E_j^* = 0$; 
\item[{\rm (iv)}]  $E_h^{\e} A E_j^{\e} = 0$; 
\item[{\rm (v)}]   $E_h^{\e} A^* E_j^{\e} = 0$.
\end{itemize}
\end{lemma}
\proof By Theorem \ref{fund:thm3a} and Lemma \ref{fund:lem4a} we have
$P \Es_h A \Es_j P^{-1} = \Ee_h \As \Ee_j$, $P^2 \Es_h A \Es_j P^{-2} = E_h \Ae E_j$, 
$P E_h \As E_j P^{-1} = \Es_h \Ae \Es_j$ and $P^2 E_h \As E_j P^{-2} = \Ee_h A \Ee_j$.
The result now follows in view of \eqref{ter:lem1i}, \eqref{ter:lem1ii} and \eqref{enakost_pr_st}. \qed

%%%%%%%%%%%%%%%%%%%%%%%%%%%%%%%%%%%%%%%%%%%%%%%%%%%%%%%%%%%%%%%%
\section{Six bases for an irreducible $\boldsymbol{T}$-module}
\label{sec:bases}

With reference to Notation \ref{blank1}, let $W$ denote an irreducible $T$-module.
We are going to show that the triple $A$, $\As$, $\Ae$ acts on $W$ as a Leonard triple.
We start with a lemma.

\begin{lemma}
\label{lem:bases}
With reference to Notation \ref{blank1}, let $W$ denote an irreducible $T$-module
with endpoint $r$ and diameter $d=D-2r$. Then the following {\rm (i)--(iii)} hold.
\begin{itemize}
\item[{\rm (i)}]   
For a nonzero $u  \in E_rW$, each of the following two sequences is a basis for $W$:
  \begin{equation}
  \label{AsA}
     \Es_r u, \Es_{r+1} u, \ldots, \Es_{r+d} u;
  \end{equation}
  \begin{equation}
  \label{AeA}
     \Ee_r u, \Ee_{r+1} u, \ldots, \Ee_{r+d} u.
  \end{equation}
\item[{\rm (ii)}]  
For a nonzero $u^* \in \Es_rW$, each of the following two sequences is a basis for $W$:
  \begin{equation}
  \label{AeAs}
     \Ee_r u^*, \Ee_{r+1} u^*, \ldots, \Ee_{r+d} u^*; 
  \end{equation}
  \begin{equation}
  \label{AAs}
     E_r u^*, E_{r+1} u^*, \ldots, E_{r+d} u^*.
  \end{equation}
\item[{\rm (iii)}] 
For a nonzero $u^\e \in \Ee_rW$, each of the following two sequences is a basis for $W$:
  \begin{equation}
  \label{AAe}
     E_r u^\e, E_{r+1} u^\e, \ldots, E_{r+d} u^\e;
  \end{equation}
  \begin{equation}
  \label{AsAe}
     \Es_r u^\e, \Es_{r+1} u^\e, \ldots, \Es_{r+d} u^\e.
  \end{equation}
\end{itemize}
\end{lemma}
\proof 
It follows from \cite[Corollaries 6.8 and 8.5]{Go} that the pair $A, \As$ acts on $W$ as a Leonard pair. 
Therefore \eqref{AsA} and \eqref{AAs} are bases for $W$ by \cite[Lemma 10.2]{ter9}. 

\smallskip \noindent
By Lemma \ref{fund:lem4a} and Corollary \ref{fund:cor4a} the sequence \eqref{AeAs} (resp. \eqref{AAe}) is 
the image under $P$ (resp. $P^2$) of the sequence \eqref{AsA}, provided $u$ is normalized so that 
$Pu=u^*$ (resp. $P^2 u = u^\e$). Similarly, the sequence \eqref{AsAe} (resp. \eqref{AeA}) is 
the image under $P$ (resp. $P^2$) of the sequence \eqref{AAs}, provided $u^*$ is normalized so that 
$Pu^*=u^\e$ (resp. $P^2 u^* = u$). Since $P$ is invertible, the sequences \eqref{AeA}, \eqref{AeAs}, 
\eqref{AAe} and \eqref{AsAe} are bases for $W$. \qed

\smallskip \noindent
The following result will be useful.

\begin{lemma}
\label{lem:bases1}
With reference to Notation \ref{blank1}, let $W$ denote an irreducible $T$-module
with endpoint $r$ and diameter $d=D-2r$. Then the following {\rm (i)--(iii)} hold.
\begin{itemize}
\item[{\rm (i)}]   For a nonzero $u \in E_rW$, 
  \begin{equation}
  \label{A_sum}
     u= \sum_{i=0}^d \Es_{r+i} u, \qquad \qquad \qquad \qquad u = \sum_{i=0}^d \Ee_{r+i} u.
  \end{equation}
\item[{\rm (ii)}]  For a nonzero $u^* \in \Es_rW$,
  \begin{equation}
  \label{As_sum}
     u^* = \sum_{i=0}^d \Ee_{r+i} u^*, \qquad \qquad \qquad \qquad u^*= \sum_{i=0}^d E_{r+i} u^*.
  \end{equation}
\item[{\rm (iii)}] For a nonzero $u^\e \in \Ee_rW$,
  \begin{equation}
  \label{Ae_sum}
     u^\e = \sum_{i=0}^d E_{r+i} u^\e, \qquad \qquad \qquad \qquad u^\e= \sum_{i=0}^d \Es_{r+i} u^\e.
  \end{equation}
\end{itemize}
\end{lemma}
\proof
(i) Evaluate $u=Iu$ using \eqref{den1} and \eqref{ter:lem2ii} to obtain the equation on the left
in \eqref{A_sum}. Evaluate $u=Iu$ using \eqref{ien2} and Lemma \ref{im_end}(i) to obtain the equation on the right in \eqref{A_sum}.

\smallskip \noindent
(ii), (iii) Similar to the proof of (i) above. \qed

\smallskip \noindent
For future use we record an idea from the proof of Lemma \ref{lem:bases}.

\begin{lemma}
\label{lem:bases2}
With reference to Notation \ref{blank1}, let $W$ denote an irreducible $T$-module
with endpoint $r$ and diameter $d=D-2r$. Then the following {\rm (i), (ii)} hold.
\begin{itemize}
\item[{\rm (i)}] 
The basis \eqref{AeAs} (resp. \eqref{AAe}, \eqref{AsA}) is the image under $P$ of the 
basis \eqref{AsA} (resp. \eqref{AeAs}, \eqref{AAe}), provided $u$ 
(resp. $u^*$, $u^\e$) is normalized so that $Pu=u^*$ (resp. $Pu^* = u^\e$, $Pu^\e = u$).
\item[{\rm (ii)}] 
The basis \eqref{AAs} (resp. \eqref{AsAe}, \eqref{AeA}) is the image under $P$ of the 
basis \eqref{AeA} (resp. \eqref{AAs}, \eqref{AsAe}), provided $u$ 
(resp. $u^*$, $u^\e$) is normalized so that $Pu=u^*$ (resp. $Pu^* = u^\e$, $Pu^\e = u$).
\end{itemize}
\end{lemma}

\begin{remark}
{\rm With reference to Notation \ref{blank1}, let $W$ denote an irreducible $T$-module
with endpoint $r$ and diameter $d=D-2r$. Note that the definition of a standard basis for $W$
\cite[Definition 6.4]{Go} is different from the 
definition of the basis \eqref{AsA}. However, it turns out that these definitions are equivalent.
Similarly, the definition of a dual standard basis for $W$ \cite[Definition 8.2]{Go} is
equivalent to the definition of the basis \eqref{AAs}. We will therefore prove all the results of 
Section \ref{sec:action} and Section \ref{sec:sp}, although some of these results were 
already proven in \cite{Go}.}
\end{remark}

%%%%%%%%%%%%%%%%%%%%%%%%%%%%%%%%%%%%%%%%%%%%%%%%%%%%%%%%%%%%%%%%

\section{The action of $\boldsymbol{A}$, $\boldsymbol{\As}$, $\boldsymbol{\Ae}$
         on the six bases}
\label{sec:action}

With reference to Notation \ref{blank1}, let $W$ denote an irreducible $T$-module with diameter 
$d$. In this section we display the matrices which represent the action of $A$, $\As$ and $\Ae$ 
on $W$ with respect to the six bases from Lemma \ref{lem:bases}. We use the following notation. 
Let $\M$ denote the $\CC$-algebra of all $(d+1) \times (d+1)$ matrices with entries in $\CC$. The 
rows and columns of matrices in $\M$ shall be indexed by $0, 1, \ldots, d$. Let 
$v_0, v_1, \ldots, v_d$ denote a basis for $W$. For $B \in \M$ and $Y \in T$ we say $B$ 
{\it represents} $Y$ {\it with respect to} $v_0, v_1, \ldots, v_d$ whenever 
$Y v_j=\sum_{i=0}^d B_{ij} v_i$ for $0 \le j \le d$. We have a comment. For an invertible $S \in T$
the following are equivalent: (i) the matrix $B$ represents $Y$ with respect to 
$v_0, v_1, \ldots, v_d$; (ii) the matrix $B$ represents $S Y S^{-1}$ with respect to 
$Sv_0, Sv_1, \ldots, Sv_d$.

\begin{theorem}
\Label{sb:thm6}
With reference to Notation \ref{blank1}, let $W$ denote an irreducible $T$-module 
with endpoint $r$ and diameter $d=D-2r$. Then the following {\rm (i)--(iii)} hold.
\begin{itemize}
\item[{\rm (i)}]  
The matrix which represents $A$ with respect to the basis {\rm \eqref{AAs}} of $W$ 
and with respect to the basis {\rm \eqref{AAe}} of $W$
is ${\rm diag}(d,d-2,d-4, \ldots, -d)$.
\item[{\rm (ii)}]  
The matrix which represents $\As$ with respect to the basis {\rm \eqref{AsA}} of $W$ 
and with respect to the basis {\rm \eqref{AsAe}} of $W$ 
is ${\rm diag}(d,d-2,d-4, \ldots, -d)$. 
\item[{\rm (iii)}]  
The matrix which represents $\Ae$ with respect to the basis {\rm \eqref{AeA}} of $W$ 
and with respect to the basis {\rm \eqref{AeAs}} of $W$
is ${\rm diag}(d,d-2,d-4, \ldots, -d)$.
\end{itemize}
\end{theorem}
\proof (i) Recall that the basis \eqref{AAs} of $W$ is of the form
$E_r u^*, E_{r+1} u^*, \ldots, E_{r+d} u^*$, where $u^*$ is a nonzero vector in $\Es_r W$. 
Similarly, the basis \eqref{AAe} of $W$ is of the form
$E_r u^\e, E_{r+1} u^\e, \ldots, E_{r+d} u^\e$, where $u^\e$ is a nonzero vector in $\Ee_rW$. 
For $0 \le i \le d$ each of $E_{r+i} u^*$, $E_{r+i} u^\e$ is an eigenvector for $A$ with 
eigenvalue $\theta_{r+i}=d-2i$. The result follows.

\smallskip \noindent
(ii), (iii) Similar to the proof of (i) above. \qed

\begin{theorem}
\Label{sb:thm7}
With reference to Notation \ref{blank1}, let $W$ denote an irreducible $T$-module 
with endpoint $r$ and diameter $d=D-2r$. Consider the following matrix in $\M$:
\begin{equation}
\label{matrika1}
\left( \begin{array}{cccccc}
	    0 & d &       &       &       & 0 \cr
            1 & 0 & d-1   &       &       &   \cr
              & 2 & \cdot & \cdot &       &   \cr
              &   & \cdot & \cdot & \cdot &   \cr
              &   &       & \cdot & \cdot & 1 \cr        
            0 &   &       &       & d     & 0 \cr
          \end{array}\right).
\end{equation}
Then the following {\rm (i)--(iii)} hold.
\begin{itemize}
\item[{\rm (i)}] 
The matrix $\eqref{matrika1}$ represents $A$ with respect to the
bases {\rm \eqref{AsA}} and {\rm \eqref{AeA}} of $W$.
\item[{\rm (ii)}]
The matrix $\eqref{matrika1}$ represents $\As$ with respect to the
bases {\rm \eqref{AeAs}} and {\rm \eqref{AAs}} of $W$.
\item[{\rm (iii)}]
The matrix $\eqref{matrika1}$ represents $\Ae$ with respect to the
bases {\rm \eqref{AAe}} and {\rm \eqref{AsAe}} of $W$.
\end{itemize}
\end{theorem}
\proof
(i)  Recall that the basis \eqref{AsA} is of the form
$\Es_r u, \Es_{r+1} u, \ldots, \Es_{r+d} u$ where $u$ is a nonzero vector in $E_rW$.
Let $B$ denote the matrix in $\M$ which represents $A$ with respect to this basis. 
We show that $B$ is equal to the matrix \eqref{matrika1}. 

For $0 \le i \le D$ we have $A \Es_i V \subseteq \Es_{i-1}V+\Es_{i+1}V$ by \eqref{ter:lem1i} 
and \eqref{odd}; therefore $A \Es_i W \subseteq \Es_{i-1}W+\Es_{i+1}W$, implying that $B$ is
tridiagonal with diagonal entries 0. Further, $(A-dI)u=0$ since $u \in E_r W$ and $\theta_r=d$.
Therefore $(B-dI)(1,1, \ldots, 1)^t = 0$ by the equation on the left in \eqref{A_sum}. 
By these comments
\begin{equation}
\label{pfeq2}
B_{0,1}=d, \qquad B_{i,i-1}+B_{i,i+1}=d \quad (1 \le i \le d-1), \qquad B_{d,d-1}=d.
\end{equation}
By \eqref{ter:lem1ii}, \eqref{enakost_pr_st} and \eqref{odd} we find
$\As E_r V \subseteq E_{r-1} V + E_{r+1} V$. By this and since $W$ has endpoint $r$ we find 
$\As E_r W \subseteq E_{r+1} W$. Therefore $\As u \in E_{r+1} W$. Now $(A-(d-2)I) \As u = 0$ since 
$\theta_{r+1}=d-2$. As $\As u = \sum_{i=0}^d \theta_{r+i} \Es_{r+i}u$ by construction, this
implies $(B-(d-2)I) (\theta_r, \theta_{r+1}, \ldots, \theta_{r+d})^t = 0$. 
Combining this with the fact that $B$ is tridiagonal with diagonal entries 0 we find
\begin{equation}
\label{pfeq3}
B_{i,i-1} \theta_{r+i-1} + B_{i,i+1} \theta_{r+i+1} = (d-2) \theta_{r+i} \qquad (1 \le i \le d-1).
\end{equation}
Combining \eqref{pfeq2} and \eqref{pfeq3} and using $\theta_{r+j} = d-2j$ for $0 \le j \le d$
we obtain $B_{i,i-1} = i$ for $1 \le i \le d$ and $B_{i,i+1}=d-i$ for $0 \le i \le d-1$. Therefore
$B$ is equal to the matrix \eqref{matrika1}.

\smallskip \noindent
Next recall that the basis \eqref{AeA} is of the form
$\Ee_r u, \Ee_{r+1} u, \ldots, \Ee_{r+d} u$ where $u$ is a nonzero vector in $E_rW$.
Let $B'$ denote the matrix in $\M$ which represents $A$ with respect to this basis. 
We show that $B'$ is equal to the matrix \eqref{matrika1}. Since the proof is similar to the 
proof that $B$ is equal to the matrix \eqref{matrika1}, we just indicate the main steps.

Similarly as in the first part of the proof, but using Lemma \ref{fund:lem4} instead of 
\eqref{ter:lem1i}, we find that $B'$ is tridiagonal with diagonal entries 0.
Since $(A-dI)u=0$ we obtain $(B'-dI)(1,1, \ldots, 1)^t = 0$ by the equation on the right in 
\eqref{A_sum}. Hence
\begin{equation}
\label{pfeq4}
B'_{0,1}=d, \qquad B'_{i,i-1}+B'_{i,i+1}=d \quad (1 \le i \le d-1), \qquad B'_{d,d-1}=d.
\end{equation}
Further, $\Ae E_r W \subseteq E_{r+1} W$ by Lemma \ref{fund:lem4} and \eqref{odd}, implying
$(A-(d-2)I) \Ae u = 0$. This gives us
$(B'-(d-2)I) (\theta_r, \theta_{r+1}, \ldots, \theta_{r+d})^t = 0$ since
$\Ae u = \sum_{i=0}^d \theta_{r+i} \Ee_{r+i}u$. Hence
\begin{equation}
\label{pfeq5}
B'_{i,i-1} \theta_{r+i-1} + B'_{i,i+1} \theta_{r+i+1} = (d-2) \theta_{r+i} 
\qquad (1 \le i \le d-1).
\end{equation}
Combining \eqref{pfeq4} and \eqref{pfeq5} we find that $B'$ is equal to the matrix 
\eqref{matrika1}.

\smallskip \noindent
(ii), (iii) 
Use Theorem \ref{fund:thm3a}, Lemma \ref{lem:bases2} and the comment above Theorem \ref{sb:thm6}. \qed

\begin{theorem}
\Label{sb:thm8}
With reference to Notation \ref{blank1}, let $W$ denote an irreducible $T$-module 
with endpoint $r$ and diameter $d=D-2r$. Consider the following matrix in $\M$:
\begin{equation}
\label{matrika2}
\im \left( \begin{array}{cccccc}
	     0 &  d &       &       &        & 0 \cr
            -1 &  0 & d-1   &       &        &   \cr
               & -2 & \cdot & \cdot &        &   \cr
               &    & \cdot & \cdot & \cdot  &   \cr
               &    &       & \cdot & \cdot  & 1 \cr        
             0 &    &       &       & -d     & 0 \cr
          \end{array}\right).
\end{equation}
Then the following {\rm (i)--(iii)} hold.
\begin{itemize}
\item[{\rm (i)}]
The matrix {\rm \eqref{matrika2}} represents $A$ with respect to the basis {\rm \eqref{AeAs}} of $W$.
\item[{\rm (ii)}]
The matrix {\rm \eqref{matrika2}} represents $\As$ with respect to the basis {\rm \eqref{AAe}} of $W$.
\item[{\rm (iii)}] 
The matrix {\rm \eqref{matrika2}} represents $\Ae$ with respect to the basis {\rm \eqref{AsA}} of $W$.
\end{itemize}
\end{theorem}
\proof (i) Let $B, B^*$ and $B^\e$ denote the matrices in $\M$ which represent $A$, $\As$ and 
$\Ae$ with respect to the basis \eqref{AeAs} of $W$. We have $B = -\im (B^* B^\e - B^\e B^*)/2$ 
by Lemma \ref{fund:lem1}(ii). Recall that $B^\e$ is equal to ${\rm diag}(d,d-2,d-4, \ldots, -d)$ by 
Theorem \ref{sb:thm6}(iii) and that $B^*$ is equal to the matrix \eqref{matrika1} by Theorem 
\ref{sb:thm7}(ii). By these comments $B$ is equal to the matrix \eqref{matrika2}. 

\smallskip \noindent
(ii), (iii) 
Use Theorem \ref{fund:thm3a}, Lemma \ref{lem:bases2} and the comment above Theorem \ref{sb:thm6}. \qed

\begin{theorem}
\Label{sb:thm9}
With reference to Notation \ref{blank1}, let $W$ denote an irreducible $T$-module 
with endpoint $r$ and diameter $d=D-2r$. Consider the following  matrix in $\M$:
\begin{equation}
\label{matrika3}
\im \left( \begin{array}{cccccc}
	    0 & -d &       &       &       &  0 \cr
            1 &  0 & 1-d   &       &       &    \cr
              &  2 & \cdot & \cdot &       &    \cr
              &    & \cdot & \cdot & \cdot &    \cr
              &    &       & \cdot & \cdot & -1 \cr        
            0 &    &       &       & d     &  0 \cr
          \end{array}\right).
\end{equation}
Then the following {\rm (i)--(iii)} hold.
\begin{itemize}
\item[{\rm (i)}]
The matrix {\rm \eqref{matrika3}} represents $A$ with respect to the basis {\rm \eqref{AsAe}} of $W$.
\item[{\rm (ii)}] 
The matrix {\rm \eqref{matrika3}} represents $\As$ with respect to the basis {\rm \eqref{AeA}} of $W$.
\item[{\rm (iii)}]
The matrix {\rm \eqref{matrika3}} represents $\Ae$ with respect to the basis {\rm \eqref{AAs}} of $W$.
\end{itemize}
\end{theorem}
\proof (i) Let $B, B^*$ and $B^\e$ denote the matrices in $\M$ which represent $A$, $\As$ and 
$\Ae$ with respect to the basis \eqref{AsAe} of $W$. We have $B = -\im (B^* B^\e - B^\e B^*)/2$ 
by Lemma \ref{fund:lem1}(ii). Recall that $B^*$ is equal to ${\rm diag}(d,d-2,d-4, \ldots, -d)$ by 
Theorem \ref{sb:thm6}(ii) and that $B^\e$ is equal to the matrix \eqref{matrika1} by Theorem 
\ref{sb:thm7}(iii). By these comments $B$ is equal to the matrix \eqref{matrika3}.

\smallskip \noindent
(ii), (iii) 
Use Theorem \ref{fund:thm3a}, Lemma \ref{lem:bases2} and the comment above Theorem \ref{sb:thm6}. \qed

\begin{corollary}
\Label{sb:cor10}
With reference to Notation \ref{blank1}, let $W$ denote an irreducible $T$-module.
Then the triple $A, \As, \Ae$ acts on $W$ as a Leonard triple.
\end{corollary}
\proof
Immediate from Theorems \ref{sb:thm6} -- \ref{sb:thm9}. \qed

%%%%%%%%%%%%%%%%%%%%%%%%%%%%%%%%%%%%%%%%%%%%%%%%%%%%%%%%%%%%%%%%%%%%%%%%%%%%%%
\section{The inner products}
\label{sec:sp}

With reference to Notation \ref{blank1}, let $W$ denote an irreducible $T$-module.
In Lemma \ref{lem:bases} we displayed six bases for $W$. Later in the paper we will find the
transition matrices between these bases. Before we get to this it is convenient to find the inner
products for the vectors in these bases.

\begin{theorem}
\Label{sp:thm1}
With reference to Notation \ref{blank1}, let $W$ denote an irreducible $T$-module 
with endpoint $r$ and diameter $d=D-2r$. Then the following {\rm (i)--(iii)} hold for $0 \le i,j \le d$.
\begin{itemize}
\item[{\rm (i)}]
For a nonzero $u \in E_rW$,
$$
  \la \Es_{r+i} u, \Es_{r+j} u \ra = \delta_{ij} {d \choose i} 2^{-d} \Vert u \Vert^2, 
  \qquad \qquad
  \la \Ee_{r+i} u, \Ee_{r+j} u \ra = \delta_{ij} {d \choose i} 2^{-d} \Vert u \Vert^2. 
$$
\item[{\rm (ii)}]
For a nonzero $u^* \in \Es_rW$,
$$
  \la \Ee_{r+i} u^*, \Ee_{r+j} u^* \ra = \delta_{ij} {d \choose i} 2^{-d} \Vert u^* \Vert^2, 
  \qquad \qquad
  \la E_{r+i} u^*, E_{r+j} u^* \ra = \delta_{ij} {d \choose i} 2^{-d} \Vert u^* \Vert^2.
$$
\item[{\rm (iii)}]
For a nonzero $u^\e \in \Ee_rW$,
$$
  \la E_{r+i} u^\e, E_{r+j} u^\e \ra = \delta_{ij} {d \choose i} 2^{-d} \Vert u^\e \Vert^2, 
  \qquad \qquad
  \la \Es_{r+i} u^\e, \Es_{r+j} u^\e \ra = \delta_{ij} {d \choose i} 2^{-d} \Vert u^\e \Vert^2.
$$
\end{itemize}
\end{theorem}
\proof
(i) Concerning the equation on the left, it holds for $i \ne j$ since $\Es_{r+i} u$ and 
$\Es_{r+j} u$ are orthogonal by \eqref{dual_direct_sum}. To verify the equation for $i=j$ we first 
claim that $\Vert \Es_{r+i}u \Vert^2 = {d \choose i} \Vert \Es_r u \Vert^2$. 
To prove the claim we assume $i \ge 1$; otherwise the result is clear.
By \eqref{cez} and since $\ov{A}^t=A$ we have 
$$
  \la A \Es_{r+i-1} u, \Es_{r+i} u \ra = \la \Es_{r+i-1} u, A \Es_{r+i} u \ra.
$$
In this equation we evaluate both sides using Theorem \ref{sb:thm7}(i) and simplify the result
using the fact that $\Es_r u, \Es_{r+1} u, \ldots, \Es_{r+d} u $ are mutually orthogonal;
we obtain
$$
  i \Vert \Es_{r+i} u \Vert^2 = (d-i+1) \Vert \Es_{r+i-1} u \Vert^2.
$$
The claim follows from this and induction on $i$. 
Next we claim that $\Vert \Es_r u \Vert^2 = 2^{-d} \Vert u \Vert^2$.
To see this, recall that  $\Es_r u, \Es_{r+1} u, \ldots, \Es_{r+d} u $ are mutually orthogonal and that
$u = \sum_{i=0}^d \Es_{r+i} u$ by the equation on the left in Lemma \ref{lem:bases1}(i). 
By these comments and the first claim,
$$
  \Vert u \Vert^2 = \sum_{i=0}^d \Vert \Es_{r+i} u \Vert^2 = 
  \Vert \Es_r u \Vert^2 \sum_{i=0}^d {d \choose i} = 2^d \Vert \Es_r u \Vert^2
$$
and the second claim is proved. Combining the two claims we get the equation on the left for $i=j$.
We have now verified the equation on the left. The proof for the equation on the right is similar,
so we just indicate the main steps. If  $i \ne j$ then $\Ee_{r+i} u$ and $\Ee_{r+j} u$ 
are orthogonal by Lemma \ref{ii:lem0}(ii). Assume now $i=j$. By \eqref{cez} we have
$\la A \Ee_{r+i-1} u, \Ee_{r+i} u \ra = \la \Ee_{r+i-1} u, A \Ee_{r+i} u \ra$. Evaluating and
simplifying this using Theorem \ref{sb:thm7}(i) and the fact that 
$\Ee_r u, \Ee_{r+1} u, \ldots, \Ee_{r+d} u $ are mutually orthogonal we obtain
$i \Vert \Ee_{r+i} u \Vert^2 = (d-i+1) \Vert \Ee_{r+i-1} u \Vert^2$. Using induction on $i$ we find
$\Vert \Ee_{r+i}u \Vert^2 = {d \choose i} \Vert \Ee_r u \Vert^2$. Using this and the equation 
on the right in Lemma \ref{lem:bases1}(i) we find $\Vert u \Vert^2 = 2^d \Vert \Ee_r u \Vert^2$.
Combining the above results we get the equation on the right.

\medskip \noindent
(ii), (iii) Use \eqref{cez}, Lemma \ref{fund:lem3a}(i), Lemma \ref{lem:bases2} and 
(i) above.  \qed

\smallskip \noindent
We have a comment.

\begin{theorem}
\Label{sp:thm3}
With reference to Notation \ref{blank1}, let $W$ denote an irreducible $T$-module 
with endpoint $r$ and diameter $d=D-2r$. Then the following {\rm (i)--(iii)} hold for $0 \le i \le d$
and nonzero vectors $u \in E_rW$, $u^* \in \Es_rW$, $u^\e \in \Ee_rW$.
\begin{itemize}
\item[{\rm (i)}] 
$E_{r+i} u^\e = \im^i (1-\im)^d \la u^\e, u^* \ra \Vert u^* \Vert^{-2} E_{r+i} u^*$.
\item[{\rm (ii)}]
$\Es_{r+i} u = \im^i (1-\im)^d \la u, u^\e \ra \Vert u^\e \Vert^{-2} \Es_{r+i} u^\e$. 
\item[{\rm (iii)}]
$\Ee_{r+i} u^* = \im^i (1-\im)^d \la u^*, u \ra \Vert u \Vert^{-2} \Ee_{r+i} u$. 
\end{itemize}
\end{theorem}
\proof 
(i) Each of $E_{r+i} u^\e$, $E_{r+i} u^*$ is a basis for $E_{r+i} W$ so there exists a nonzero
$\lambda_i \in \CC$ such that $E_{r+i} u^\e = \lambda_i E_{r+i} u^*$. We first claim that 
$\lambda_i = \im^i \lambda_0$. To prove the claim assume $i \ge 1$; otherwise the result is
clear. Note that $\As E_{r+i} u^\e = \lambda_i \As E_{r+i} u^*$. Evaluating both sides of 
this equation using Theorem \ref{sb:thm7}(ii) and Theorem \ref{sb:thm8}(ii) and then comparing the 
results we find that $\lambda_i = \im \lambda_{i-1}$. The claim follows from this and induction
on $i$. Next we claim that $\lambda_0 = (1-\im)^d \la u^\e, u^* \ra \Vert u^* \Vert^{-2}$.
To see this recall that $u^\e = \sum_{i=0}^d E_{r+i} u^\e$ and $u^* = \sum_{i=0}^d E_{r+i} u^*$
by Lemma \ref{lem:bases1}(ii),(iii). By this, \eqref{de}, Theorem \ref{sp:thm1}(ii) and the first
claim we find 
$$
  \la u^\e, u^* \ra = \sum_{i=0}^d \la E_{r+i} u^\e, E_{r+i} u^* \ra = 
  \lambda_0 \sum_{i=0}^d \im^i \Vert E_{r+i} u^* \Vert^2 = 
  {\lambda_0 \Vert u^* \Vert^2 \over 2^d} \sum_{i=0}^d {d \choose i} \im^i =
  {\lambda_0 \Vert u^* \Vert^2 (1+\im)^d \over 2^d}
$$
and the second claim follows. Combining the two claims we obtain the desired result. 

\smallskip \noindent
(ii), (iii) Use \eqref{cez}, Lemma \ref{fund:lem3a}(i), Lemma \ref{lem:bases2} and 
(i) above. \qed

\begin{corollary}
\Label{sp:cor3}
With reference to Notation \ref{blank1}, let $W$ denote an irreducible $T$-module 
with endpoint $r$. Let $u$, $u^*$, $u^\e$ denote nonzero vectors 
in $E_rW$, $\Es_rW$, $\Ee_rW$, respectively. Then each of $\la u^\e, u^* \ra$, $\la u, u^\e \ra$,
$\la u^*, u \ra$ is nonzero.
\end{corollary}
\proof
The vectors $E_r u^\e$, $\Es_r u$ and $\Ee_r u^*$ are nonzero by Lemma \ref{lem:bases}. 
Combining this with Theorem \ref{sp:thm3} we get the result. \qed

\begin{theorem}
\Label{sp:thm2}
With reference to Notation \ref{blank1}, let $W$ denote an irreducible $T$-module 
with endpoint $r$ and diameter $d=D-2r$. Then the following {\rm (i)--(iii)} hold for $0 \le i,j \le d$
and nonzero vectors $u \in E_rW$, $u^* \in \Es_rW$, $u^\e \in \Ee_rW$.
\begin{itemize}
\item[{\rm (i)}] 
$\la E_{r+i} u^\e , E_{r+j} u^* \ra = 
 \delta_{ij} \: \im^i {d \choose i} (1+\im)^{-d} \la u^\e, u^* \ra$.
\item[{\rm (ii)}] 
$\la \Es_{r+i} u, \Es_{r+j} u^\e \ra = 
\delta_{ij} \: \im^i {d \choose i} (1+\im)^{-d} \la u, u^\e \ra$.
\item[{\rm (iii)}] 
$\la \Ee_{r+i} u^*, \Ee_{r+j} u \ra = 
\delta_{ij} \: \im^i {d \choose i} (1+\im)^{-d} \la u^*, u \ra$.
\end{itemize}
\end{theorem}
\proof
Immediate from Theorem \ref{sp:thm1} and Theorem \ref{sp:thm3}. \qed

\smallskip \noindent
Before proceeding we recall a definition. For an integer $n \ge 0$ and $a \in \CC$ we define
$$
  (a)_n = a(a+1)(a+2) \cdots (a+n-1).
$$
We interpret $(a)_0=1$. For integers $0 \le i,j \le d$ we define
\begin{equation}
\label{hyper}
_2F_1 \Big( {-i,-j \atop -d} ; 2 \Big) = \sum_{n=0}^d {(-i)_n (-j)_n \over (-d)_n \, n!} \, 2^n.
\end{equation}
The sum \eqref{hyper} is an example of a {\em hypergeometric series} \cite[Section 2.1]{AAR}. 
We will use the fact that 
\begin{equation}
\label{enakost}
_2F_1 \Big( {-i,-j \atop -d};2 \Big)={d-2j \over d-i+1} \: _2F_1 \Big( {1-i,-j \atop -d};2 \Big) -
{i-1 \over d-i+1} \; _2F_1 \Big( {2-i,-j \atop -d};2 \Big)
\end{equation}
provided $i \ge 2$. Line \eqref{enakost} follows from \cite[Equation 1.10.3]{KS} and since each 
side of \eqref{hyper} is equal to $K_i(j;1/2,d)$ \cite[Definition 1.10.1]{KS} where the  
$K_n(x;p,N)$ are the Krawtchouk polynomials.

\begin{theorem}
\label{sp:thm4}
With reference to Notation \ref{blank1}, let $W$ denote an irreducible $T$-module 
with endpoint $r$ and diameter $d=D-2r$. Then the following {\rm (i)--(iii)} hold for $0 \le i,j \le d$
and nonzero vectors $u \in E_rW$, $u^* \in \Es_rW$, $u^\e \in \Ee_rW$.
\begin{itemize}
\item[{\rm (i)}]
$$
  \la E_{r+i} u^*, \Es_{r+j} u \ra = 
  2^{-d} \la u^*, u \ra {d \choose i} {d \choose j} \: _2F_1 \Big( {-i,-j \atop -d};2 \Big).
$$
\item[{\rm (ii)}]
$$
  \la \Es_{r+i} u^\e, \Ee_{r+j} u^* \ra = 
  2^{-d} \la u^\e, u^* \ra {d \choose i} {d \choose j} \: _2F_1 \Big( {-i,-j \atop -d};2 \Big).
$$
\item[{\rm (iii)}]
$$
  \la \Ee_{r+i} u, E_{r+j} u^\e \ra = 
  2^{-d} \la u, u^\e \ra {d \choose i} {d \choose j} \: _2F_1 \Big( {-i,-j \atop -d};2 \Big).
$$
\end{itemize}
\end{theorem}
\proof 
(i) We first claim that
\begin{equation}
\label{wv*a}
  \la E_{r+i} u^*, \Es_{r+j} u \ra =
  \la E_r u^*, \Es_r u \ra {d \choose i} {d \choose j} \: _2F_1 \Big( {-i,-j \atop -d};2 \Big).
\end{equation}
We will follow the approach of Go \cite[Theorem 9.1]{Go} and prove the claim using induction
on $i+j$. It is clear that \eqref{wv*a} holds for $i=j=0$. Assume now $d \ge 1$; otherwise we are done. 
To show that \eqref{wv*a} holds for $(i,j)=(0,1)$ observe by \eqref{cez}, 
Theorem \ref{sb:thm6}(i) and Theorem \ref{sb:thm7}(i) that
\begin{equation}
\label{eq10}
\la E_r u^*, \Es_{r+1} u \ra = \la E_r u^*, A \Es_r u \ra = 
\la A E_r u^*, \Es_r u \ra =  d \la E_r u^*, \Es_r u \ra.
\end{equation}
Similarly \eqref{wv*a} holds for $(i,j)=(1,0)$ since
$$
  \la E_{r+1} u^*, \Es_r u \ra = \la \As E_r u^*, \Es_r u \ra = 
  \la E_r u^*, \As \Es_r u \ra = d \la E_r u^*, \Es_r u \ra.
$$
To show that \eqref{wv*a} holds for $(i,j)=(1,1)$, observe by \eqref{cez}, Theorem \ref{sb:thm6}(ii)
and Theorem \ref{sb:thm7}(ii) that
$$
  \la E_{r+1} u^*, \Es_{r+1} u \ra = \la \As E_r u^*, \Es_{r+1} u \ra =
  \la E_r u^*, \As \Es_{r+1} u \ra = (d-2) \la E_r u^*, \Es_{r+1} u \ra,
$$
and that the last expression is equal to $d (d-2) \la E_r u^*,\Es_r u \ra$ by \eqref{eq10}.

\smallskip \noindent
For the rest of this proof assume $d \ge 2$ and $i \ge 2$ or $j \ge 2$; otherwise we are done. 
We first assume $i \ge 2$. By Theorem \ref{sb:thm7}(ii),
$$
  E_{r+i} u^* = {1 \over i} \As E_{r+i-1} u^* - {d-i+2 \over i} E_{r+i-2} u^*.
$$
Using this, \eqref{cez} and Theorem \ref{sb:thm6}(ii) we obtain
\begin{equation}
\label{enacba}
  \la E_{r+i} u^*, \Es_{r+j} u \ra = 
  {d-2j \over i} \la E_{r+i-1} u^*, \Es_{r+j} u \ra - 
  {d-i+2 \over i} \la E_{r+i-2} u^*,\Es_{r+j} u \ra.
\end{equation}
By the induction hypothesis the right-hand side of \eqref{enacba} is equal to
${d \choose i} {d \choose j} \la E_r u^*, \Es_r u \ra$ times
$$
  {d-2j \over d-i+1} \: _2F_1 \Big( {1-i,-j \atop -d};2 \Big) -
  {i-1 \over d-i+1} \: _2F_1 \Big( {2-i,-j \atop -d};2 \Big).
$$
Evaluating the above expression using \eqref{enakost} we obtain \eqref{wv*a}.
Now assume $j \ge 2$. By Theorem \ref{sb:thm7}(i),
$$
  \Es_{r+j} u = {1 \over j} A \Es_{r+j-1} u - {d-j+2 \over j} \Es_{r+j-2} u.
$$
Using this, \eqref{cez} and Theorem \ref{sb:thm6}(i) we obtain
\begin{equation}
\label{enacba1}
  \la E_{r+i} u^*, \Es_{r+j} u \ra = 
  {d-2i \over j} \la E_{r+i} u^*, \Es_{r+j-1} u \ra - 
  {d-j+2 \over j} \la E_{r+i} u^*,\Es_{r+j-2} u \ra.
\end{equation}
By the induction hypothesis the right-hand side of \eqref{enacba1} is equal to
${d \choose i} {d \choose j} \la E_r u^*, \Es_r u \ra$ times
$$
  {d-2i \over d-j+1} \: _2F_1 \Big( {-i,1-j \atop -d};2 \Big) -
  {j-1 \over d-j+1} \: _2F_1 \Big( {-i,2-j \atop -d};2 \Big).
$$
Evaluating the above expression using \eqref{enakost} we obtain \eqref{wv*a} 
and our first claim is proved. 
Next we claim that $\la E_r u^*, \Es_r u \ra = 2^{-d} \la u^*,u \ra$. To see this, observe
by \eqref{cez}, \eqref{den3}, \eqref{den2}, Lemma \ref{lem:bases1}(ii) and \eqref{wv*a} that
$$
  \la u^*, u \ra = \la \Es_r u^*, u \ra = \la u^*, \Es_r u \ra = 
  \sum_{i=0}^d \la E_{r+i} u^*, \Es_r u \ra =
  \la E_r u^*, \Es_r u \ra \sum_{i=0}^d {d \choose i} = \la E_r u^*, \Es_r u \ra \: 2^d
$$
and the second claim follows. Combining the two claims we get the desired result.

\smallskip \noindent
(ii), (iii) Use \eqref{cez}, Lemma \ref{fund:lem3a}(i), Lemma \ref{lem:bases2} and 
(i) above. \qed

\begin{theorem}
\label{sp:thm5}
With reference to Notation \ref{blank1}, let $W$ denote an irreducible $T$-module 
with endpoint $r$ and diameter $d=D-2r$. Then the following {\rm (i)--(iii)} hold for $0 \le i,j \le d$
and nonzero vectors $u \in E_rW$, $u^* \in \Es_rW$, $u^\e \in \Ee_rW$.
\begin{itemize}
\item[{\rm (i)}]
$$
  \la E_{r+i} u^*, \Es_{r+j} u^\e \ra = 
  \im^j \, 2^{-d} \la u^*, u^\e \ra {d \choose i} {d \choose j} 
  \: _2F_1 \Big( {-i,-j \atop -d};2 \Big).
$$
\item[{\rm (ii)}]
$$
  \la \Es_{r+i} u^\e, \Ee_{r+j} u \ra = 
  \im^j \, 2^{-d} \la u^\e, u \ra {d \choose i} {d \choose j} 
  \: _2F_1 \Big( {-i,-j \atop -d};2 \Big).
$$
\item[{\rm (iii)}]
$$
  \la \Ee_{r+i} u , E_{r+j} u^* \ra = 
  \im^j \, 2^{-d} \la u, u^* \ra  {d \choose i} {d \choose j} 
  \: _2F_1 \Big( {-i,-j \atop -d};2 \Big).
$$
\end{itemize}
\end{theorem}
\proof
(i) 
By Theorem \ref{sp:thm3}(ii) and \eqref{wv*a},
\begin{equation}
\label{seenaena}
  \la E_{r+i} u^*, \Es_{r+j} u^\e \ra = 
  \im^j \la E_r u^*, \Es_r u^\e \ra {d \choose i} {d \choose j} 
  \: _2F_1 \Big( {-i,-j \atop -d};2 \Big).
\end{equation}
We claim that $\la E_r u^*, \Es_r u^\e \ra = 2^{-d} \la u^*, u^\e \ra$. 
To see this, observe
by \eqref{cez}, \eqref{den3}, \eqref{den2}, Lemma \ref{lem:bases1}(ii) and \eqref{seenaena} that
\begin{small}
$$
  \la u^*, u^\e \ra = \la \Es_r u^*, u^\e \ra = \la u^*, \Es_r u^\e \ra = 
  \sum_{i=0}^d \la E_{r+i} u^*, \Es_r u^\e \ra =
  \la E_r u^*, \Es_r u^\e \ra \sum_{i=0}^d {d \choose i} = \la E_r u^*, \Es_r u^\e \ra 2^d
$$
\end{small}
and the claim follows.
Combining \eqref{seenaena} with the claim we get the desired result.

\smallskip \noindent
(ii), (iii) Use \eqref{cez}, Lemma \ref{fund:lem3a}(i), Lemma \ref{lem:bases2} and 
(i) above. \qed 

\begin{theorem}
\label{sp:thm7}
With reference to Notation \ref{blank1}, let $W$ denote an irreducible $T$-module 
with endpoint $r$ and diameter $d=D-2r$. Then the following {\rm (i)--(iii)} hold for $0 \le i,j \le d$
and nonzero vectors $u \in E_rW$, $u^* \in \Es_rW$, $u^\e \in \Ee_rW$.
\begin{itemize}
\item[{\rm (i)}]
$$
  \la E_{r+i} u^\e, \Es_{r+j} u \ra = 
  \im^i \, 2^{-d} \la u^\e, u \ra {d \choose i} {d \choose j} 
  \: _2F_1 \Big( {-i,-j \atop -d};2 \Big).
$$
\item[{\rm (ii)}]
$$
  \la \Es_{r+i} u, \Ee_{r+j} u^* \ra = 
  \im^i \, 2^{-d} \la u, u^* \ra {d \choose i} {d \choose j} 
  \: _2F_1 \Big( {-i,-j \atop -d};2 \Big).
$$
\item[{\rm (iii)}]
$$
  \la \Ee_{r+i} u^*, E_{r+j} u^\e \ra = 
  \im^i \, 2^{-d} \la u^*, u^\e \ra {d \choose i} {d \choose j} 
  \: _2F_1 \Big( {-i,-j \atop -d};2 \Big).
$$
\end{itemize}
\end{theorem}
\proof (i) By Theorem \ref{sp:thm3}(i),(ii) and \eqref{seenaena},
\begin{equation}
\label{seenaenaena}
  \la E_{r+i} u^\e, \Es_{r+j} u \ra = 
    \im^i \la E_r u^\e, \Es_r u \ra {d \choose i} {d \choose j} 
  \: _2F_1 \Big( {-i,-j \atop -d};2 \Big).
\end{equation}
We claim that $\la E_r u^\e, \Es_r u \ra = 2^{-d} \la u^\e, u \ra$. 
To see this, observe
by \eqref{cez}, \eqref{pi3}, \eqref{pi4}, Lemma \ref{lem:bases1}(i) and \eqref{seenaenaena} that
$$
  \la u^\e, u \ra = \la u^\e, E_r u \ra =  \la E_r u^\e, u \ra = 
  \sum_{j=0}^d \la E_r u^\e, \Es_{r+j} u \ra =
  \la E_r u^\e, \Es_r u \ra \sum_{j=0}^d {d \choose j} = \la E_r u^\e, \Es_r u \ra \: 2^d
$$
and the claim follows. Combining \eqref{seenaenaena} with the claim we get the desired result.
   
\smallskip \noindent
(ii), (iii) Use \eqref{cez}, Lemma \ref{fund:lem3a}(i), Lemma \ref{lem:bases2} and 
(i) above. \qed

\begin{theorem}
\label{sp:thm6}
With reference to Notation \ref{blank1}, let $W$ denote an irreducible $T$-module 
with endpoint $r$ and diameter $d=D-2r$. Then the following {\rm (i)--(iii)} hold for $0 \le i,j \le d$
and nonzero vectors $u \in E_rW$, $u^* \in \Es_rW$, $u^\e \in \Ee_rW$.
\begin{itemize}
\item[{\rm (i)}]
$$
  \la E_{r+i} u^*, \Ee_{r+j} u^* \ra = 
  \im^{-i-j} (2-2\im)^{-d} \Vert u^* \Vert^2 {d \choose i} {d \choose j} 
  \: _2F_1 \Big( {-i,-j \atop -d};2 \Big).
$$
\item[{\rm (ii)}]
$$
  \la \Es_{r+i} u^\e, E_{r+j} u^\e \ra = 
  \im^{-i-j} (2-2\im)^{-d} \Vert u^\e \Vert^2 {d \choose i} {d \choose j} 
  \: _2F_1 \Big( {-i,-j \atop -d};2 \Big).
$$
\item[{\rm (iii)}]
$$
  \la \Ee_{r+i} u, \Es_{r+j} u \ra = 
  \im^{-i-j} (2-2\im)^{-d} \Vert u \Vert^2 {d \choose i} {d \choose j} 
  \: _2F_1 \Big( {-i,-j \atop -d};2 \Big).
$$
\end{itemize}
\end{theorem}
\proof (i) 
By Theorem \ref{sp:thm3}(i),
$$
  \la E_{r+i} u^*, \Ee_{r+j} u^* \ra = 
  {\Vert u^* \Vert^2 \over \im^i(1-\im)^d \la u^\e, u^* \ra} \la E_{r+i} u^\e, \Ee_{r+j} u^* \ra.
$$
The result now follows from Theorem \ref{sp:thm7}(iii).

\smallskip \noindent
(ii), (iii) Use \eqref{cez}, Lemma \ref{fund:lem3a}(i), Lemma \ref{lem:bases2} and 
(i) above. \qed

%%%%%%%%%%%%%%%%%%%%%%%%%%%%%%%%%%%%%%%%%%%%%%%%%%%%%%%%%%%%%%%%%%%%%%%%%%%%%%
\section{The inner products between $\mathbf{u}$, $\mathbf{u^*}$ and $\mathbf{u^\e}$}
\label{sec:spu}

With reference to Notation \ref{blank1}, let $W$ denote an irreducible $T$-module 
with endpoint $r$, and pick nonzero vectors $u \in E_rW$, $u^* \in \Es_rW$,
$u^\e \in \Ee_rW$. In this section we display some equations
involving the inner products between $u, u^*, u^\e$.

\begin{theorem}
\label{spu:thm1}
With reference to Notation \ref{blank1}, let $W$ denote an irreducible $T$-module 
with endpoint $r$ and diameter $d=D-2r$. Then the following {\rm (i)--(iii)} hold for nonzero vectors 
$u \in E_rW$, $u^* \in \Es_rW$, $u^\e \in \Ee_rW$.
\begin{itemize}
\item[{\rm (i)}]
$\Vert u \Vert^2 = (1+\im)^d \la u, u^* \ra \la u^\e,u \ra \la u^\e,u^* \ra^{-1}$.
\item[{\rm (ii)}]
$\Vert u^* \Vert^2 = (1+\im)^d \la u^*, u^\e \ra \la u,u^* \ra \la u,u^\e \ra^{-1}$.
\item[{\rm (iii)}]
$\Vert u^\e \Vert^2 = (1+\im)^d \la u^\e, u \ra \la u^*,u^\e \ra \la u^*,u \ra^{-1}$.
\end{itemize}
\end{theorem}
\proof (i) 
Using Theorem \ref{sp:thm3}(iii) and Theorem \ref{sp:thm4}(ii) we obtain
$$
  \la \Es_r u^\e, \Ee_r u \ra = 
  {\Vert u \Vert^2 \, \la u^\e, u^* \ra \over 2^d \, (1+\im)^d \, \la u,u^* \ra}.
$$
Comparing the above value for $\la \Es_r u^\e, \Ee_r u \ra$ with the value given in 
Theorem \ref{sp:thm5}(ii) we obtain the desired result. 

\smallskip \noindent
(ii), (iii) Similar to the proof of (i) above. \qed

\begin{corollary}
\label{spu:cor1a}
With reference to Notation \ref{blank1}, let $W$ denote an irreducible $T$-module 
with endpoint $r$ and diameter $d=D-2r$. Pick nonzero vectors 
$u \in E_rW$, $u^* \in \Es_rW$, $u^\e \in \Ee_rW$. Then the scalar
\begin{equation}
\label{stevilo}
  \la u, u^* \ra \la u^*, u^\e \ra \la u^\e, u \ra (1+\im)^d
\end{equation}
is real and positive.
\end{corollary}
\proof By Theorem \ref{spu:thm1}(i) the scalar \eqref{stevilo} is equal to 
$\Vert u \Vert^2 \la u^*, u^\e \ra \la u^\e, u^* \ra$. By construction $\Vert u \Vert^2$ is real
and positive. Also $\la u^*, u^\e \ra \la u^\e, u^* \ra$ is real and positive since
$\la u^\e, u^* \ra = \ov{\la u^*, u^\e \ra}$ by construction and $\la u^\e, u^* \ra \ne 0$
by Corollary \ref{sp:cor3}. The result follows. \qed

\smallskip \noindent
With reference to Theorem \ref{spu:thm1}, the inner products $\la u,u^* \ra$, $\la u^*,u^\e \ra$ and 
$\la u^\e, u \ra$ are independent in the following sense.

\begin{lemma}
\label{spu:lem2}
Let $a,b,c \in \CC$ be such that $a b c (1+\im)^d$ is a positive real number. Then 
there exist nonzero vectors $u \in E_r W$, $u^* \in \Es_r W$, $u^\e \in \Ee_r W$ such 
that $\la u,u^* \ra=a$, $\la u^*,u^\e \ra=b$, $\la u^\e, u \ra=c$. 
\end{lemma}
\proof
Note that each of $a,b,c$ is nonzero.
Pick arbitrary nonzero vectors $u_1 \in E_r W$, $u_1^* \in \Es_r W$, $u_1^\e \in \Ee_r W$
and define $\delta=a b c (1+\im)^d \Vert c u_1^* \Vert^{-2}$. 
Also define $\lambda = a \, \delta^{-1/2} \la u_1,u_1^* \ra^{-1}$, 
$\lambda^* = \delta^{1/2}$, $\lambda^\e = \ov{b} \, \delta^{-1/2} \la u_1^\e,u_1^* \ra^{-1}$. 
These scalars are well-defined by Corollary \ref{sp:cor3}. 
We now define $u= \lambda u_1$, $u^* = \lambda^* u_1^*$, $u^\e = \lambda^\e u_1^\e$ and
routinely obtain $\la u,u^* \ra=a$, $\la u^*,u^\e \ra=b$, $\la u^\e, u \ra=c$ using
Theorem \ref{spu:thm1}. \qed

%%%%%%%%%%%%%%%%%%%%%%%%%%%%%%%%%%%%%%%%%%%%%%%%%%%%%%%%%%%%%%%%%%%%%%%%%%%%%%
\section{The transition matrices}
\label{sec:tra}

In this section we display the transition matrices between the bases introduced in Section
\ref{sec:bases}. We start with a comment. With reference to Notation \ref{blank1}, let $W$ denote 
an irreducible $T$-module with endpoint $r$ and diameter $d=D-2r$. Let $u_0, \ldots, u_d$ and 
$v_0, \ldots, v_d$ denote bases for $W$. 
By the {\em transition matrix from} $u_0, \ldots, u_d$ {\em to} $v_0, \ldots, v_d$ 
we mean the matrix $C \in \M$ which satisfies 
$$
  v_j = \sum_{i=0}^d C_{ij} u_i \qquad (0 \le j \le d).
$$
We recall a few properties of transition matrices. Let $C$ denote the transition matrix from 
$u_0, \ldots, u_d$ to $v_0, \ldots, v_d$. Then $C$ is invertible and $C^{-1}$ is 
the transition matrix from $v_0, \ldots, v_d$ to $u_0, \ldots, u_d$.
If $u_0, \ldots, u_d$ are mutually orthogonal 
%and $v_0, \ldots, v_d$ are mutually orthogonal
then the entries of $C$ are given by
\begin{equation}
\label{trmat}
  C_{ij} = {\la v_j, u_i \ra \over \Vert u_i \Vert^2} \qquad (0 \le i,j \le d).
\end{equation}
In order to display the transition matrices in a compact form we abbreviate 
$$  
  \Phi_{ij}= {d \choose j} \: _2F_1 \Big( {-i,-j \atop -d};2 \Big)
$$
and
$$
  D_1 = {\rm diag}(\im^0, \im^1, \ldots, \im^d), \qquad \qquad 
  D_2 = {\rm diag}(\im^{-0}, \im^{-1}, \ldots, \im^{-d}).
$$

\begin{theorem}
\label{tm:tm}
With reference to Notation \ref{blank1}, let $W$ denote an irreducible $T$-module
with endpoint $r$ and diameter $d=D-2r$. Pick nonzero vectors
$u \in E_rW$, $u^* \in \Es_rW$, $u^\e \in \Ee_rW$.
Then the transition matrices between the bases \eqref{AsA} -- \eqref{AsAe} are given in the
tables below.

%\begin{footnotesize}
\begin{center}
\begin{tabular}{|c|c|c|c|}
% prva vrstica
\hline
 & 
to the basis \eqref{AsA} & 
to the basis \eqref{AeA} & 
to the basis \eqref{AeAs} \\

% druga vrstica
\hline
\begin{tabular}{cc}
from the \\ basis \eqref{AsA}
\end{tabular}                  & 
${\rm diag}(1,1, \ldots, 1)$     & 
$ {1 \over (1-\im)^d} \big[ \im^{-i-j} \; \Phi_{ij} \big]_{i,j=0}^d $ &
${\la u^*, u \ra \over \Vert u \Vert^2} \big[ \im^{-i} \; \Phi_{ij} \big]_{i,j=0}^d$ \\ 

% tretja vrstica
\hline
\begin{tabular}{cc}
from the \\ basis \eqref{AeA}
\end{tabular}                  & 
${1 \over (1+\im)^d} \big[ \im^{i+j} \; \Phi_{ij} \big]_{i,j=0}^d$ & 
${\rm diag}(1, 1, \ldots, 1)$ & 
${(1-\im)^d \la u^*, u \ra \over \Vert u \Vert^2} \, D_1$ \\
  
% cetrta vrstica
\hline
\begin{tabular}{cc}
from the \\ basis \eqref{AeAs}
\end{tabular}                  & 
${\la u, u^* \ra \over \Vert u^* \Vert^2} 
  \big[ \im^j \; \Phi_{ij} \big]_{i,j=0}^d$  & 
${(1+\im)^d \la u, u^* \ra \over \Vert u^* \Vert^2} \, D_2$ & 
${\rm diag}(1, 1, \ldots, 1)$ \\

%peta vrstica
\hline
\begin{tabular}{cc}
from the \\ basis \eqref{AAs}
\end{tabular}                  & 
${\la u, u^* \ra \over \Vert u^* \Vert^2} 
  \big[ \Phi_{ij} \big]_{i,j=0}^d$ & 
$  {\la u, u^* \ra \over \Vert u^* \Vert^2} 
  \big[ \im^i \; \Phi_{ij} \big]_{i,j=0}^d$ &
${1 \over (1+\im)^d}
  \big[ \im^{i+j} \; \Phi_{ij} \big]_{i,j=0}^d $\\
  
%sesta vrstica
\hline
\begin{tabular}{cc}
from the \\ basis \eqref{AAe}
\end{tabular}                  &
$ {\la u, u^\e \ra \over \Vert u^\e \Vert^2} 
  \big[ \im^{-i} \; \Phi_{ij} \big]_{i,j=0}^d$ &
$ {\la u, u^\e \ra \over \Vert u^\e \Vert^2} 
  \big[ \Phi_{ij} \big]_{i,j=0}^d$ & 
${\la u^*, u^\e \ra \over \Vert u^\e \Vert^2} 
  \big[ \im^j \; \Phi_{ij} \big]_{i,j=0}^d $ \\
  
%sedma vrstica
\hline
\begin{tabular}{cc}
from the \\ basis \eqref{AsAe}
\end{tabular}                  &
${(1-\im)^d \la u, u^\e \ra \over \Vert u^\e \Vert^2} \, D_1$ &
${\la u, u^\e \ra \over \Vert u^\e \Vert^2} 
  \Big[ \im^{-j} \; \Phi_{ij} \big]_{i,j=0}^d$ &
${\la u^*, u^\e \ra \over \Vert u^\e \Vert^2} 
  \big[ \Phi_{ij} \big]_{i,j=0}^d$ \\
\hline
\end{tabular}

\bigskip 

\begin{tabular}{|c|c|c|c|}
% prva vrstica
\hline
 & 
to the basis \eqref{AAs} & 
to the basis \eqref{AAe} & 
to the basis \eqref{AsAe} \\

% druga vrstica
\hline
\begin{tabular}{cc}
from the \\ basis \eqref{AsA}
\end{tabular}                  & 
${\la u^*, u \ra \over \Vert u \Vert^2} 
  \big[ \Phi_{ij} \big]_{i,j=0}^d$     & 
${\la u^\e, u \ra \over \Vert u \Vert^2} 
  \big[ \im^j \; \Phi_{ij} \big]_{i,j=0}^d $ &
${(1+\im)^d \la u^\e, u \ra \over \Vert u \Vert^2} \, D_2$ \\ 

%tretja vrstica
\hline
\begin{tabular}{cc}
from the \\ basis \eqref{AeA}
\end{tabular}                  &
${\la u^*, u \ra \over \Vert u \Vert^2} 
  \big[ \im^{-j} \; \Phi_{ij} \big]_{i,j=0}^d$ &
${\la u^\e, u \ra \over \Vert u \Vert^2} 
  \big[ \Phi_{ij} \big]_{i,j=0}^d$ &
$ {\la u^\e, u \ra \over \Vert u \Vert^2} 
  \big[ \im^i \; \Phi_{ij} \big]_{i,j=0}^d $ \\

%cetrta vrstica
\hline
\begin{tabular}{cc}
from the \\ basis \eqref{AeAs}
\end{tabular}                  &
${1 \over (1-\im)^d} 
  \big[ \im^{-i-j} \; \Phi_{ij} \big]_{i,j=0}^d $ &
${\la u^\e, u^* \ra \over \Vert u^* \Vert^2} 
  \big[ \im^{-i} \; \Phi_{ij} \big]_{i,j=0}^d$ &
${\la u^\e, u^* \ra \over \Vert u^* \Vert^2} 
  \big[ \Phi_{ij} \big]_{i,j=0}^d $ \\

%peta vrstica
\hline
\begin{tabular}{cc}
from the \\ basis \eqref{AAs}
\end{tabular}                  &
${\rm diag}(1, 1, \ldots, 1)$ &
${(1-\im)^d \la u^\e, u^* \ra \over \Vert u^* \Vert^2} \, D_1$ &
${\la u^\e, u^* \ra \over \Vert u^* \Vert^2} 
  \big[ \im^{-j} \; \Phi_{ij} \big]_{i,j=0}^d$ \\

%sesta vrstica
\hline
\begin{tabular}{cc}
from the \\ basis \eqref{AAe}
\end{tabular}                  &
${(1+\im)^d \la u^*, u^\e \ra \over \Vert u^\e \Vert^2} \, D_2$ &
${\rm diag}(1, 1, \ldots, 1)$ &
${1 \over (1-\im)^d}
  \big[ \im^{-i-j} \; \Phi_{ij} \big]_{i,j=0}^d$ \\

%sedma vrstica
\hline
\begin{tabular}{cc}
from the \\ basis \eqref{AsAe}
\end{tabular}                  &
${\la u^*, u^\e \ra \over \Vert u^\e \Vert^2} 
  \big[ \im^i \; \Phi_{ij} \big]_{i,j=0}^d $ &
${1 \over (1+\im)^d}
  \big[ \im^{i+j} \; \Phi_{ij} \big]_{i,j=0}^d$ &
${\rm diag}(1, 1, \ldots, 1)$ \\
\hline
\end{tabular}
\end{center}
%\end{footnotesize}
\end{theorem}
\proof
Combine Theorem \ref{sp:thm1}, Theorems \ref{sp:thm2} -- \ref{sp:thm6} and \eqref{trmat}. \qed

\bigskip \noindent
{\bf \large Acknowledgement:} The author would like to thank Paul Terwilliger
for proposing the problem and for his careful reading of the earlier versions 
of this manuscript.

%%%%%%%%%%%%%%%%%%%%%%%%%%%%%%%%%%%%%%%%%%%%%%%%%%%%%%%%%%%%%%%%%%%%%%%%%%%%%%%%%%%%

\end{document}